\documentclass[3p,authoryear,fleqn]{elsarticle}

\usepackage{ecrc}
\usepackage{algorithm}
\usepackage{algpseudocode}
\usepackage{amsmath}
\usepackage{graphicx}
\graphicspath{ {./images/} }
\usepackage{longtable}
\usepackage{array}
\usepackage{booktabs,adjustbox}
\usepackage{epstopdf}
\usepackage{amsmath}
\usepackage[utf8]{inputenc}
\usepackage{amssymb}
\usepackage{color}
\usepackage{algorithm,algpseudocode}
\usepackage{caption}
\usepackage{array}
\usepackage{float}
\usepackage{comment}
\usepackage{lscape}
\usepackage{booktabs}
\allowdisplaybreaks
\usepackage[colorinlistoftodos,prependcaption,textsize=tiny]{todonotes}
\usepackage[title]{appendix}
\usepackage{hyperref}
\usepackage[official]{eurosym}
\usepackage{multirow}
\usepackage{layouts}
\usepackage{placeins}
\usepackage{afterpage}
\usepackage{tabularx}
\newcolumntype{L}{>{\raggedright\arraybackslash}X}
\newcolumntype{s}{>{\hsize=.15\hsize}X}

\volume{00}

\firstpage{1}

\journalname{Journal of Cleaner Production}

\runauth{Belyak et al.}

\usepackage{amssymb}
\usepackage[figuresright]{rotating}

\begin{document}

\begin{frontmatter}

%% Title, authors and addresses

%% use the tnoteref command within \title for footnotes;
%% use the tnotetext command for the associated footnote;
%% use the fnref command within \author or \address for footnotes;
%% use the fntext command for the associated footnote;
%% use the corref command within \author for corresponding author footnotes;
%% use the cortext command for the associated footnote;
%% use the ead command for the email address,
%% and the form \ead[url] for the home page:
%%
%% \title{Title\tnoteref{label1}
%% \tnotetext[label1]{}
%% \author{Name\corref{cor1}\fnref{label2}
%% \ead{email address}
%% \ead[url]{home page}
%% \fntext[label2]{}
%% \cortext[cor1]{}
%% \address{Address\fnref{label3}
%% \fntext[label3]{}

% \dochead{}
%% Use \dochead if there is an article header, e.g. \dochead{Short communication}
%% \dochead can also be used to include a conference title, if directed by the editors
%% e.g. \dochead{17th International Conference on Dynamical Processes in Excited States of Solids}

%TC:ignore

% Main title of the paper
\title{Renewable Energy Expansion under Taxes and Subsidies: A Transmission Operator's Perspective}

%% use optional labels to link authors explicitly to addresses:
%% \author[label1,label2]{<author name>}
%% \address[label1]{<address>}
%% \address[label2]{<address>}

\address[Aalto]{Aalto University, Espoo, 02150, Finland}
\address[UM]{University of Maryland, College Park, Maryland, 20742, United States}
\address[NTNU]{Norwegian University of Science and Technology, Trondheim, NO-7491, Norway}
\address[IIASA]{International Institute for Applied Systems Analysis, Laxenburg, A-2361, Austria}

\author[Aalto]{Nikita Belyak}
\author[UM,Aalto,NTNU]{Steven A. Gabriel}
\author[IIASA]{Nikolay Khabarov}
\author[Aalto]{Fabricio Oliveira}

\begin{abstract}
This paper investigates the role of a transmission system operator within a carbon footprint reduction strategy incorporating carbon taxes and renewable energy generation subsidies in the decentralised energy market. This is achieved via an optimisation bi-level model in which a welfare-maximizing transmission system operator makes investments in transmission lines at the upper level while considering power market dynamics at the lower level. To account for the deregulated energy market structure, this paper assumes that the generation companies at the lower level make capacity investments as price-takers in perfect competition. Considering alternative transmission infrastructure expansion budgets, carbon emission taxes and monetary incentives for renewable energy generation capacity expansion, the impact of alternative compositions of these factors is analysed against three output factors: the share of renewable energy in the generation mix, total generation amount, and social welfare. The proposed modelling assessment is applied to an illustrative three-node instance and a case study considering a simplified representation of the energy system of the Nordic and Baltic countries. The results highlight that, under certain circumstances, renewable energy generation subsidies may lead to an increase of renewable energy in the generation mix followed by a simultaneous fall in the total generation amount. Nevertheless, when applied together, these three measures demonstrated a positive impact on all output factors within Nordics' and Baltics' energy systems. The experiments additionally suggest that considering the high value of the carbon tax does not have an impact on the output factors while the composition of high values of renewable energy generation subsidies and budget for transmission infrastructure expansion has the strongest effect. %The findings of this study can serve as an additional reference for policymakers seeking to simultaneously increase the proportion of renewable energy generation and total generation amount to meet increasing demand. 
%Nevertheless, one should take into account that due to the limitations of modern computational software, the proposed modelling approach involved several simplifications in both the detailing of the energy system and the reduction of input parameter space. These shortcomings should be addressed in further research. 
\end{abstract}

\begin{keyword}
%% keywords here, in the form: keyword \sep keyword
Transmission planning \sep Renewable generation \sep Energy modelling \sep Carbon tax \sep renewable generation subsidies

%% PACS codes here, in the form: \PACS code \sep code

%% MSC codes here, in the form: \MSC code \sep code
%% or \MSC[2008] code \sep code (2000 is the default)

\end{keyword}

%TC:endignore

\end{frontmatter}

\pagebreak

\section{Introduction}

The rise in the world's population and associated increasing environmental footprint has underscored the pressing need for global sustainable development. In 2015, the United Nations identified 17 pillars for sustainable development. \citet{taghvaee2023sustainability} suggest that one of the most essential pillars that have an impact on others is global openness and international agreements. A parallel perspective is echoed by \citet{nasrollahi2020environmental}, who attempted to bring the attention of policymakers to international agreements, technological development and energy efficiency emphasising their pivotal role in various sustainable development pillars. This pivotal role of energy efficiency is further underscored by \citet{taghvaee2022comparing}, who highlighted its substantial impact on social welfare, the economy and the environment. 

The prospect of ensuring a high level of energy efficiency closely aligns with the need for action to address the alarming effects of climate change. As a response, many European Union (EU) countries have undergone fundamental transformations of their power generation sector to expand the share of low-carbon renewable energy \citep{wolf_european_2021, noauthor_climate_nodate, STEFFEN2018280, Cambridge_climate}. Despite a considerable amount of renewable power that has been installed in the past decade, the vast majority of the EU member countries are far from meeting these levels of variable renewable energy sources (VREs) individually \citep{agora_european_2021}. Consequentially, a significant amount of renewable energy capacity will be built in the short- and medium-term, requiring large-scale investments.   

However, in liberalised electricity markets, such as those found in EU countries, the UK, and North America, renewable energy revenues are insufficient to provide an adequate return to VRE capacity for private investors \citep{HAAR2020111483}. Additionally, the increase of the VRE share among the energy sources may jeopardise sufficient power quality (i.e., requirements for uninterrupted power supply and stable conditions of voltage and current), energy systems stability, power balance and efficient transmission and distribution of power \citep{SINSEL20202271}. However, existing transmission systems design is not capable of coping with significant levels of renewable penetration \citep{7752917}. Consequentially, renewable-driven expansion of generation requires new approaches for transmission network planning.

The surveys on existing transmission expansion planning literature \cite{HEMMATI2013312, LUMBRERAS201619, 7583779} suggest the decisions regarding the structure of the power market, the level of detail on the operation of the system and the solution method for the problem to be amongst the key factors defining the distinct approaches. Mathematical optimisation has been extensively applied as a solution method primarily in an academic context as it eradicates the risk of suboptimality of the solution \cite{LUMBRERAS201619}. However, one should take into account the trade-off between the computational feasibility of solving the problem to optimality and its scale, which in turn is augmented as the modelling detail level and the size of the network modelled increase.  

Regarding the representation of deregulation in the power sector, i.e., decoupling transmission and generation expansion decisions, one can pinpoint two generalised strategies in the literature. The first spans investigations aimed at developing an optimal transmission network expansion strategy that would account for various possible developments of the generation infrastructure. Examples of such strategy can be found in \citep{SUN2018546, MORTAZ201935}. Nonetheless, while the burden to formulate exhaustive uncertainty sets appears to be challenging on its own, this modelling strategy prevents generation companies (GenCos) from being dynamic market players capable of making reactive decisions regarding generation levels and capacity expansion. 

Another strategy attempts to develop efficient modelling tools to consider the planning of the transmission and generation infrastructure expansion in a coordinated manner. For example, this coordinated modelling approach has been considered in \citep{7752917, 8945251, ZHANG2020105944}. For the modelling assessment proposed in this paper, we consider a decentralised planning strategy to ensure the representation of the reactive position (i.e., acting as price-takers) of GenCos.  

Therefore, this paper aims to study the influence of the transmission system operator (TSO)’s decisions in composition with existing renewables-driven policies, such as carbon taxes and VRE-associated incentives, on the generation mix while assuming the GenCos to behave as price-takers in perfect-competition market settings. It's important to emphasize that these three policies can be considered mechanisms of good governance, which has been shown to play a crucial role in ensuring the efficient management of natural resources \citep{tatar2024good}. Furthermore, this study concentrates on assessing the effectiveness of various policy combinations on the generation mix and social welfare in a decentralised energy market involving generation companies of different scales. To evaluate the impact of the aforementioned renewables-driven policies in the generation mix, applied both separately and together, we propose and utilise a bi-level optimisation model. 

The proposed model assumes the TSO to take a leading position and anticipate the generation capacity investment decisions influenced by its transmission system expansion. This assumption leads to the bi-level structure of the proposed model. Such a modelling approach is widely used in energy market planning. As an example, \citet{ZHANG201684} exploited a bi-level scheme to consider integrated generation-transmission expansion at the upper level and modified unit-commitment model with demand response at the lower level. \citet{VIRASJOKI2020104716} considered a bi-level structure when formulating the model for optimal energy storage capacity sizing and use planning. In this paper, we reformulate the model proposed in \citep{VIRASJOKI2020104716} to consider welfare maximising TSO at the upper level, making decisions in the transmission lines instead of energy storage. An analogous strategy has been considered by \citet{SIDDIQUI2019208} during the investigation of the indirect influence of the TSO’s decisions as a part of an emissions mitigation strategy aligned with different levels of carbon charges in a deregulated industry. Aimed at the analytical implications, their paper neglects VRE and demand-associated uncertainty, as well as the heterogeneity of the GenCos, while assuming unlimited generation capacity. These shortcomings are addressed in the current paper by means of introducing VRE intermittency and allowing GenCos to invest in diversified power generation technologies. Furthermore, we account for various investment budget portfolios for TSO and GenCos to investigate how GenCos' investment capital availability influences the total VRE share in the optimal generation mix. 

This paper evaluates the efficiency of renewables-driven policies by considering a combined transmission and generation investment modelling assessment that involves 1) welfare-maximising TSO at the upper level, 2) heterogeneous GenCos acting under perfect competition at the lower level, 3) uncertainty associated with VRE availability, 4) carbon tax and 5) investment budget constraints to closely represent realistic power system expansion problem. Different levels of the carbon tax and incentives supporting VRE capacity deployment, as well as budget levels allocated for transmission and generation infrastructure expansion, are considered input parameters. The structure of the proposed model allows us to investigate whether TSO's decisions may efficiently foster renewable generation capacity expansion on its own or if it would require other instruments, such as carbon tax and incentives supporting VRE capacity investment. Moreover, we also evaluate the effectiveness of the latter two individually and in combination with the other policies. 

The proposed modelling approach allows one to identify optimal strategies for an illustrative 3-node energy system instance and a more realistic energy system replicating the features of Nordic and Baltic countries. In the numerical experiments, a sensitivity analysis is conducted by evaluating the changes in optimal total welfare as well as VRE share in the optimal generation mix and the optimal amount of energy generated. Therefore, the outcome of this research fills in some important knowledge gaps regarding the role of the TSO in a decarbonisation strategy.

The paper is structured as follows. In Section \ref{sec: problem formulation}, we present the optimisation problem formulations considering centralised and decentralised power market structure. Section \ref{sec: mppdc} presents the primal-dual approach utilised to transform bi-level problems given in Section \ref{sec: problem formulation} into a tractable single-level version. In Sections \ref{sec: illustrative_example} and \ref{sec: nordic_case}, we evaluate the efficiency of the proposed modelling framework by applying it to two case studies. Section \ref{sec: conclusion} presents conclusions. 

\section{Problem formulation}
\label{sec: problem formulation}
In the following sections, we present two alternative formulations of the transmission infrastructure expansion planning problem, considering the centralised and decentralised energy market structure, where the latter is modelled as perfect competition. The notation utilised is described in detail in Appendix \ref{appendix: list of nomenclature}. 
\subsection{Centralised planning}
\label{section: centralised_market}

We begin with a centralised market structure modelled as a single-level problem that is used as a benchmark in this paper. With such a setting, the central planner aims at maximising social welfare by means of making both optimal transmission ($l^+_{n,m}$) and conventional ($g^{e+}_{n,i}$) and renewable ($g^{r+}_{n,i}$) generation capacity investments decisions. It is important to highlight that pre-existing generation and transmission capacity prior to the beginning of the planning horizon is also accounted for in the modelling setup. The VRE availability is modelled via the consideration of different scenarios $s \in S$ (or operation modes, as sometimes referred to in the literature) and assumptions on the percentage of the total VRE capacity available at each of the nodes considering each of the scenarios and time periods $t \in T$. This value is referred to as the availability factor and denoted as $A^r_{s,t,n}$ for each of the VRE types. Each of the scenarios is weighted with probability $P_s$ such that $\sum_{s \in S} P_s = 1$. The parameter $D^e$ defines the carbon tax (\euro/MWh) for the conventional generation of type $e \in E$. Likewise, the VRE incentive ${INC}_n$ is defined for each of the nodes $n \in N$ as the \% value that is deducted from the VRE generation capacity investment costs represented by the parameter $I^r_{n,i}$. Therefore, if one, for example, assumed the VRE incentive ${INC}_n$ to be 10\% the GenCos will have to pay 90\% of the investment costs associated with the VRE capacity expansion.

With the aforementioned assumptions in mind, the objective function for the centralised planning model can be formulated as follows:
\begin{align}
   \nonumber \max &  \sum_{ n \in N} \Bigg( \sum_{ t \in T}  \sum_{s \in S}  P_s \left[ D^{int}_{s,t,n}q_{s,t,n} - \frac{1}{2}D^{slp}_{s,t,n} \frac{q_{s,t,n}^2}{T_t}  - \sum_{i \in I } 
    \sum_{e \in E} \left(C^{e}_{n,i} + D^e \right) g^{e}_{s,t,n,i}  \right] \\ \nonumber
   % & + \sum_{i \in I}  \pi  \left( CO^2_i - \sum_{n \in N }\sum_{t \in T }\sum_{s \in S }\sum_{e \in E }P_s D^e g^e_{s,t,n,i}   \right) \\ \nonumber
   &  -\sum_{i \in I }\Bigg( \sum_{r \in R}\left[ M^r_{n,i} \left(\overline{G}^r_{n,i} + \overline{g}^{r+}_{n,i}\right) + (1 - \frac{{INC}_n}{100}) I^r_{n,i}\overline{g}^{r+}_{n,i} \right]   \\ \nonumber
   & + \sum_{e \in E} \left[ M^{e}_{n,i}   \left( \overline{G}^{e}_{n,i} + \overline{g}^{e +}_{n,i} \right) + I^{e}_{n,i} \overline{g}^{e+}_{n,i} \right] \Bigg) \\   &  -  \sum_{m \in N} \left[ M^l_{n,m}\frac{1}{2}\left(\overline{L}_{n,m}+l^+_{n,m}\right) + \frac{1}{2} I^l_{n,m} l^+_{n,m}\right] 
   \Bigg).  \label{slm: objective}
   \end{align}
   
   The objective function \eqref{slm: objective} maximises the total profit for all the GenCos and TSO. The variable $q_{s,t,n}$ represents the demand at node $n$ during time period $t$ and considering scenario $s$. Therefore, $\forall s\in S, \ t \in T, \ n \in N$ the term $ \left(D^{int}_{s,t,n} - \frac{1}{2}D^{slp}_{s,t,n} \frac{q_{s,t,n}}{T_t} \right) q_{s,t,n} $ represents the revenue, where the fist multiplier correspond to the energy price at node $n$ during time period $t$ and considering scenario $s$.
   
   The conventional and renewable generation capacity expansion decisions ($g^{e+}_{n,i}$ and $g^{r+}_{n,i}$, respectively) along with the transmission lines capacity expansion decisions ($l^+_{n,m}$) are modelled as continuous variables (in MW) restricted by the existing investment budget limits ($K^{g+}_{i}$ and $K^{L+}$, respectively) and non-negativity conditions as follows: 
   \begin{align} 
    & l^+_{n,m} - l^+_{m,n} = 0 \quad \forall n \in N, m \in N \label{slm: lines_invest_primal_feas_1}\\  
      & \sum_{n \in N} \sum_{m \in N}I^l_{n,m}l^+_{n,m} - K^{L+} \le 0 \label{slm: lines_invest_budget_limit}\\
          & \sum_{n \in N}\sum_{r \in R} I^r_{n,i} \overline{g}^{r+}_{n,i} + \sum_{n \in N}\sum_{e \in E} I^e_{n,i} \overline{g}^{e+}_{n,i} - K^{g+}_{i} \le 0 \quad \forall i \in I \label{slm: gen_exp_budget}\\
        & \overline{g}^{e+}_{n,i} \ge 0 \quad \forall e \in E, n \in N, i \in I \label{slm: ge+_non-negativity}\\
   & \overline{g}^{r+}_{n,i} \ge 0 \quad \forall r \in R, n \in N, i \in I \label{slm: gr+_non-negativity} \\
   & l^+_{n,m} \ge 0 \quad \forall n \in N,m \in N \label{slm: l+_non-negativity},
   \end{align}
   where \eqref{slm: lines_invest_primal_feas_1} guarantees the equivalence of the capacity expansion decisions for the transmission lines $n \rightarrow m$ and $m \rightarrow n $, respectively, $\forall n,m \in N$. It is important to highlight that the continuous nature of the variables representing capacity expansion ($g^{e+}_{n, i}$, $g^{r+}_{n, i}$ and $l^+_{n,m}$) is a simplification needed to limit computational requirements and make feasible the computational experiments performed. 
   
   The conventional and renewable generation levels ($g^e_{s,t,n,i}$ and $g^r_{s,t,n,i}$, respectively) are decided considering the following constraints
   \begin{align} \nonumber
      &  q_{s,t,n} - \sum_{i \in I} \left[\sum_{e \in E}g^{e}_{s,t,n,i} + \sum_{r \in R} g^{r}_{s,t,n,i}\right] \\ 
    & +  \sum_{m \in N: m > n} f_{s,t, n, m}  -  \sum_{m \in N: m < n} f_{s, t, m, n} = 0  \quad  \forall  s \in S, t \in T, n \in N  \label{slm: balance} \\
    & g^{r}_{s,t,n,i} - T_t A^r_{s,t,n} \left(\overline{G}^{r}_{n,i} + \overline{g}^{r +}_{n,i} \right) \le 0 \quad \forall s \in S, t \in T, n \in N, i \in I, r \in R   \label{slm: VRES_gen_bounds}\\
    & g^{e}_{s,t,n,i} - T_t \left(\overline{G}^{e}_{n,i} + \overline{g}^{e+}_{n,i}\right) \le 0 \quad \forall s \in S, t \in T,  n \in N, i \in I,  e \in E   \label{slm: conv_gen_bounds}\\
     & g^e_{s,t,n,i} - g^e_{s,t-1,n,i} - T_t R^{up,e}_{n,i} \left( \overline{G}^e_{n,i} + \overline{g}^{e+}_{n.i}\right) \le 0 \quad \forall s \in S, t\in T, n \in N, i \in I, e \in E,  \label{slm: ramp-up}\\
     & g^e_{s,t-1,n,i} - g^e_{s,t,n,i} - T_t R^{down,e}_{n,i}\left( \overline{G}^e_{n,i} + \overline{g}^{e+}_{n.i}\right) \le 0 \quad \forall  s \in S, t\in T, n \in N, i \in I, e \in E,  \label{slm: ramp-down} \\
     & q_{s,t,n} \ge 0 \quad  \forall   s \in S, t \in T, n \in N\\
     & g^{e}_{s,t,n,i} \ge 0 \quad \forall s \in S, t \in T, n \in N, i \in I, e \in E \label{slm: ge_non-negativity}\\
    & g^{r}_{s,t,n,i} \ge 0 \quad \forall s \in S,  t \in T, n \in N, i \in I, r \in R,    \quad  \label{slm: gr_non-negativity}
\end{align}
where constraint \eqref{slm: balance} ensures the balance between power consumption, generation and transmission, while constraints \eqref{slm: VRES_gen_bounds} and \eqref{slm: conv_gen_bounds} define the bounds for the VRE and conventional generation, respectively. Inequalities \eqref{slm: ramp-up} and \eqref{slm: ramp-down} represent maximum ramping levels for conventional generation.
 
The terms $f_{s,t,n,m}$ represent the power flows from node $n$ to node $m$ ( $f_{s,t,n,m}>0$ if node $n$ supplies energy to node $m$, and $f_{s,t,n,m}=0$ otherwise) and are constrained by the following set of constraints 
\begin{align}
    & f_{s,t,n,m} - T_t \left(\overline{L}_{n,m} + l^+_{n,m} \right) \le 0\quad  \forall s \in S,  t \in T,  n \in N, m \in N \label{slm: flow_bounds_1}\\
    & - f_{s,t,n,m} - T_t \left(\overline{L}_{n,m} + l^+_{n,m} \right) \le 0\quad  \forall s \in S, t \in T, n \in N, m \in N   \label{slm: flow_bounds_2}\\
    & f_{s,t,n,m} = 0 \quad \forall s \in S, t \in T, n \in N, m \in N : m \le n \label{slm: flow_primal_feas}, 
\end{align}
where the equality \eqref{slm: flow_primal_feas} excludes from optimisation the terms $f_{s,t,n,m}$ where $m \le n$, constraints \eqref{slm: flow_bounds_1} and \eqref{slm: flow_bounds_2} define the transmission flows bounds. It is worth highlighting that we opted for simplifying the transmission operation system details for the sake of computational tractability. Therefore, power flow dynamics only takes into account Kirchhoff's first law stating that the sum of currents entering each node should be zero \cite{4909474} rather than a DC or AC flow model, which ultimately circumvents the computational burden incurred from the need to consider additional control variables \cite{LUMBRERAS201619, PADIYAR198817}.

In the following sections, we separately formulate lower and upper portions of the proposed bi-level optimisation problem representing a decentralised power market structure, as opposed to the centralised structure presented in Section \ref{section: centralised_market}. 
\subsection{Lower level: power market operations}

The lower level of the proposed bi-level formulation depicts the power market operations under the assumption of having been given the TSO's decisions regarding the capacity of the transmission lines. The modelling setting considers the possibility for each GenCo to own and operate some generation capacity of various technologies at each node $n \in N$ simultaneously and sell the power generated to the market. Hence, each GenCo invests in the additional generation capacity to maximise its revenue. 

We assume the market to experience perfect competition, implying that GenCos act as price-takes. Following \citep[Section 3.4.2]{gabriel_complementarity_2013} and considering linear inverse demand functions, we can formulate a single optimisation problem objective function as follows. 
\begin{align}
   \nonumber \max &  \sum_{ n \in N} \Bigg( \sum_{ t \in T}  \sum_{s \in S}  P_s \bigg[ D^{int}_{s,t,n}q_{s,t,n} - \frac{1}{2}D^{slp}_{s,t,n}\frac{q_{s,t,n}^2}{T_t}  - \sum_{i \in I } 
    \sum_{e \in E} \left( C^{e}_{n,i} + D^e \right)g^{e}_{s,t,n,i} \bigg]  \\ \nonumber  
    %& + \sum_{i \in I}  \pi  \left( CO^2_i - \sum_{n \in N }\sum_{t \in T }\sum_{s \in S }\sum_{e \in E }P_s D^e g^e_{s,t,n,i}   \right) \\ \nonumber
   &  -\sum_{i \in I }\Bigg( \sum_{r \in R}\left[ M^r_{n,i} \left(\overline{G}^r_{n,i} + \overline{g}^{r+}_{n,i}\right) + (1 - \frac{{INC}_n}{100})I^r_{n,i}\overline{g}^{r+}_{n,i} \right]   \\ 
   & + \sum_{e \in E} \left[ M^{e}_{n,i}   \left( \overline{G}^{e}_{n,i} + \overline{g}^{e +}_{n,i} \right) + I^{e}_{n,i} \overline{g}^{e+}_{n,i} \right] \Bigg) \Bigg) \label{llm: objective} \\ \nonumber
   \nonumber \text{s.t.: } & \eqref{slm: gen_exp_budget} - \eqref{slm: gr+_non-negativity}, \eqref{slm: balance} - \eqref{slm: flow_primal_feas}.
\end{align}

\subsection{Upper-level: transmission infrastructure planning}

The formulation of the upper level of the proposed bi-level problem defining the optimal transmission expansion strategy planned by the TSO is given as follows. Considering a proactive TSO's position, it attempts to maximise social welfare through transmission infrastructure expansion, taking into account a given investment capability while anticipating GenCos' reactive generation capacity expansion and operation decisions at the lower level. Hence, the corresponding optimisation problem can be formulated as 
\begin{align}
   \nonumber \max_{l^+_{n,m}} &  \sum_{ n \in N} \Bigg( \sum_{ t \in T}  \sum_{s \in S}  P_s \Bigg[ D^{int}_{s,t,n}q_{s,t,n} - \frac{1}{2}D^{slp}_{s,t,n}q_{s,t,n}^2 & \\ \nonumber
   & - \sum_{i \in I } 
    \sum_{e \in E} \left( \left[C^{e}_{n,i} + D^{e} \right]g^{e}_{s,t,n,i} \right) \Bigg] \\\nonumber
   &  -\sum_{i \in I }\Bigg( \sum_{r \in R}\left[ M^r_{n,i} \left(\overline{G}^r_{n,i} + \overline{g}^{r+}_{n,i}\right) + (1 - \frac{{INC}_n}{100})I^r_{n,i}\overline{g}^{r+}_{n,i} \right]   \\ \nonumber
   & + \sum_{e \in E} \left[ M^{e}_{n,i}   \left( \overline{G}^{e}_{n,i} + \overline{g}^{e +}_{n,i} \right) + I^{e}_{n,i} \overline{g}^{e+}_{n,i} \right] \Bigg) \\ \  &  -  \sum_{m \in N} \left[ M^l_{n,m}\frac{1}{2}\left(\overline{L}_{n,m}+l^+_{n,m}\right) + \frac{1}{2}I^l_{n,m}l^+_{n,m}\right]\Bigg)   \label{ulm: objective} \\ \nonumber
   \text{s.t.: }  &  \eqref{slm: lines_invest_primal_feas_1} -  \eqref{slm: lines_invest_budget_limit}, \eqref{slm: l+_non-negativity}  \\ \nonumber &  \text{and} \ 
   g^e_{s,t,n,i}, q_{s,t,n}, \overline{g}^{e+}_{n,i}, \overline{g}^{r+}_{n,i} \in \ \text{arg} \max \left\{\eqref{llm: objective} \text{ s.t.: } \eqref{slm: gen_exp_budget} - \eqref{slm: gr+_non-negativity}, \eqref{slm: balance} - \eqref{slm: flow_primal_feas} \right\}.
\end{align}

\section{Single-level representation of the bi-level problem}
\label{sec: mppdc}

To solve the proposed bi-level model
\begin{align}
\nonumber & \max \big\{ \eqref{ulm: objective}, \ \text{s. t.: } \ \eqref{slm: lines_invest_primal_feas_1} -  \eqref{slm: lines_invest_budget_limit}, \eqref{slm: l+_non-negativity}, \  \text{and} \\ & 
   g^e_{s,t,n,i}, \ g^r_{s,t,n,i}, \ q_{s,t,n}, \ \overline{g}^{e+}_{n,i}, \ \overline{g}^{r+}_{n,i} \in \ \text{arg} \max \left\{\eqref{llm: objective}, \ \text{s. t.: } \eqref{slm: gen_exp_budget} - \eqref{slm: gr+_non-negativity}, \eqref{slm: balance} - \eqref{slm: flow_primal_feas} \right\} \big\} 
   \label{bi_level_formulation}
 \end{align}
 we employ a single-level reformulation. In contrast to \citep{VIRASJOKI2020104716} who developed a single-level reformulation of an analogous problem based on a mathematical program with primal and dual constraints approach, we relied on a formulation that yields a mathematical program with equilibrium constraints (MPEC) \citep{gabriel_complementarity_2013}, which involves constraints that represent equilibrium conditions. Considering that the lower-level optimization problem is convex and all the constraints are affine functions, one can replace it with its corresponding Karush–Kuhn–Tucker (KKT) conditions. The KKT conditions, in this context, represent the lower-level problem's optimality conditions. Then, the MPEC comprises the upper-level problem objective function \eqref{ulm: objective} and respective constraints \eqref{slm: lines_invest_primal_feas_1}, \eqref{slm: lines_invest_budget_limit} and \eqref{slm: l+_non-negativity} with constraints representing the KKT conditions for the lower-level problem. Notice that the KKT conditions for the lower-level problem contain primal feasibility constraints $\eqref{slm: gen_exp_budget} - \eqref{slm: gr+_non-negativity}$, $\eqref{slm: balance} - \eqref{slm: flow_primal_feas}$ and the following set of constraints:
\begin{align}
      \nonumber &  \nabla_{g^e_{s,t,n,i}} \eqref{llm: objective} + \theta_{s,t,n} \nabla_{g^e_{s,t,n,i}} \eqref{slm: balance} + \beta^{e}_{s,t,n,i}  \nabla_{g^e_{s,t,n,i}}  \eqref{slm: conv_gen_bounds} \\ & + \beta^{up, e}_{s,t,n,i}  \nabla_{g^e_{s,t,n,i}}  \eqref{slm: ramp-up} + \beta^{down, e}_{s,t,n,i}  \nabla_{g^e_{s,t,n,i}}  \eqref{slm: ramp-down}  - \lambda^e_{s,t,n,i}= 0, \  \forall \ s \in S, t \in T, n \in N, i \in I, e \in E\\ 
             &  \nabla_{g^r_{s,t,n,i}} \eqref{llm: objective} + \theta_{s,t,n} \nabla_{g^r_{s,t,n,i}} \eqref{slm: balance} + \beta^{r}_{s,t,n,i}  \nabla_{g^r_{s,t,n,i}}  \eqref{slm: VRES_gen_bounds} - \lambda^r_{s,t,n,i} = 0, \  \forall \ s \in S, t \in T, n \in N, i \in I, r \in R  \\ 
                   &  \nabla_{q_{s,t,n}} \eqref{llm: objective} + \theta_{s,t,n} \nabla_{q_{s,t,n}} \eqref{slm: balance} = 0, \  \forall \ s \in S, t \in T, n \in N,  \\ 
                        \nonumber &  \nabla_{\overline{g}^{e+}_{n,i}} \eqref{llm: objective} +  \beta^{g+}_{i}\nabla_{\overline{g}^{e+}_{n,i}}\eqref{slm: gen_exp_budget} + \sum_{s \in S} \sum_{t\in T}  \beta^{e}_{s,t,n,i}  \nabla_{\overline{g}^{e+}_{n,i}}  \eqref{slm: conv_gen_bounds} \\  & + \sum_{s \in S} \sum_{t \in T}  \beta^{up, e}_{s,t,n,i}  \nabla_{\overline{g}^{e+}_{n,i}}  \eqref{slm: ramp-up}  + \sum_{s \in S} \sum_{t\in T} \beta^{down, e}_{s,t,n,i}\nabla_{\overline{g}^{e+}_{n,i}} - \lambda^{e+}_{n,i}= 0, \ \forall \ n \in N, I \in I, e \in E  \\ 
                              &  \nabla_{\overline{g}^{r+}_{n,i}} \eqref{llm: objective} +  \sum_{i } \beta^{g+}_{i}\nabla_{\overline{g}^{r+}_{n,i}}\eqref{slm: gen_exp_budget}  + \sum_{s \in S} \sum_{t\in T}  \beta^{r}_{s,t,n,i}  \nabla_{\overline{g}^{r+}_{n,i}}  \eqref{slm: VRES_gen_bounds}- \lambda^{r+}_{n,i} = 0, \ \forall \ n \in N, I \in I, r \in R  \\ 
                              \nonumber &  \theta_{s,t,n} \nabla_{
                              {f}^{s,t,n,m}_{n,i}} \eqref{slm: balance}  +  \beta^{f_1}_{s,t,n,m} \nabla_{{f}^{s,t,n,m}_{n,i}}   \eqref{slm: flow_bounds_1}  \\ & + \beta^{f_2}_{s,t,n,m} \nabla_{{f}^{s,t,n,m}_{n,i}} \eqref{slm: flow_bounds_2} + \lambda^{f}_{s,t,n,m} \nabla_{{f}^{s,t,n,m}_{n,i}} \eqref{slm: flow_primal_feas} = 0, \ \forall \ s \in S, t \in T, n \in N, m \in N : m \le n \\ 
      & 0 \le \beta^{g+}_{i} \perp \eqref{slm: gen_exp_budget} \le 0,  \ \forall \ i \in I\\ 
      & 0 \le \beta^{r}_{s,t,n,i} \perp \eqref{slm: VRES_gen_bounds} \le 0,  \ \forall \ s \in S, t \in T, n \in N, i \in I, r \in R \\
      & 0 \le \beta^{e}_{s,t,n,i} \perp \eqref{slm: conv_gen_bounds} \le 0,  \ \forall \ s \in S, t \in T, n \in N, i \in I, e \in E \\
      & 0 \le \beta^{up, e}_{s,t,n,i} \perp \eqref{slm: ramp-up} \le 0,  \ \forall \ s \in S, t \in T, n \in N, i \in I, e \in E \\
      & 0 \le \beta^{down, e}_{s,t,n,i} \perp \eqref{slm: ramp-down} \le 0,  \ \forall \ \ s \in S, t \in T, n \in N, i \in I, e \in E \\
      &  0 \le  \beta^{f_1}_{s,t,n,m} \perp \eqref{slm: flow_bounds_1} \le 0,  \ \forall \ s \in S, t \in T, n \in N, m \in N  \\
        &  0 \le  \beta^{f_2}_{s,t,n,m} \perp \eqref{slm: flow_bounds_2} \le 0,  \ \forall \ s \in S, t \in T, n \in N, m \in N  \\
        & 0 \le \lambda^{e+}_{n,i} \perp -\eqref{slm: ge+_non-negativity} \le 0,  \ \forall \ n \in N, i \in I, e \in E \\ 
            & 0 \le \lambda^{r+}_{n,i} \perp -\eqref{slm: gr+_non-negativity} \le 0,  \ \forall \ n \in N, i \in I, r \in R \\ 
             & 0 \le \lambda^{e}_{s,t,n,i} \perp -\eqref{slm: ge_non-negativity} \le 0,  \ \forall \ n \in N, i \in I, e \in E \\ 
             & 0 \le \lambda^{r}_{s,t,n,i} \perp -\eqref{slm: gr_non-negativity} \le 0,  \ \forall \ n \in N, i \in I, r \in R \\ 
              & \text{and }\eqref{slm: flow_primal_feas}, \nonumber
\end{align}
%
%where $\Omega = \left\{ g^e_{s,t,n,i}, \ g^r_{s,t,n,i}, \ q_{s,t,n}, \ \overline{g}^{e+}_{n,i}, \ \overline{g}^{r+}_{n,i} \right\}$. 
where $\nabla_x(\cdot)$ represents a vector whose components are partial deriatives of $~\cdot~$ with respect to $x$.

\section{Numerical experiments: Illustrative example}
%\label{sec: experiments}
\label{sec: illustrative_example}

In this section, the proposed modelling assessment is conducted by applying it to the first example. The resulting optimal decision policies and their impact on a perfectly competitive market structure are compared with the benchmark model, i.e., the centralised planner optimal strategy. 
All the models are implemented using the Julia (version 1.3.1) language \citep{Bazanson2017Julia} utilising the BilevelJuMP.jl package \citep{Joaquim_bilevel_2022} and solved using the commercial solver Gurobi (version
9.0.0) \citep{gurobi}. The source code and data generated for the illustrative three-node and the Nordics case studies are openly available at GitHub repositories \citep{belyak_ill_2022, belyak_nor_2022}.

%\subsection{Illustrative example: three-node system}
%\label{sec: illustrative_example}

\begin{figure}[h]
    \centering
  \includegraphics[width=0.7\linewidth]{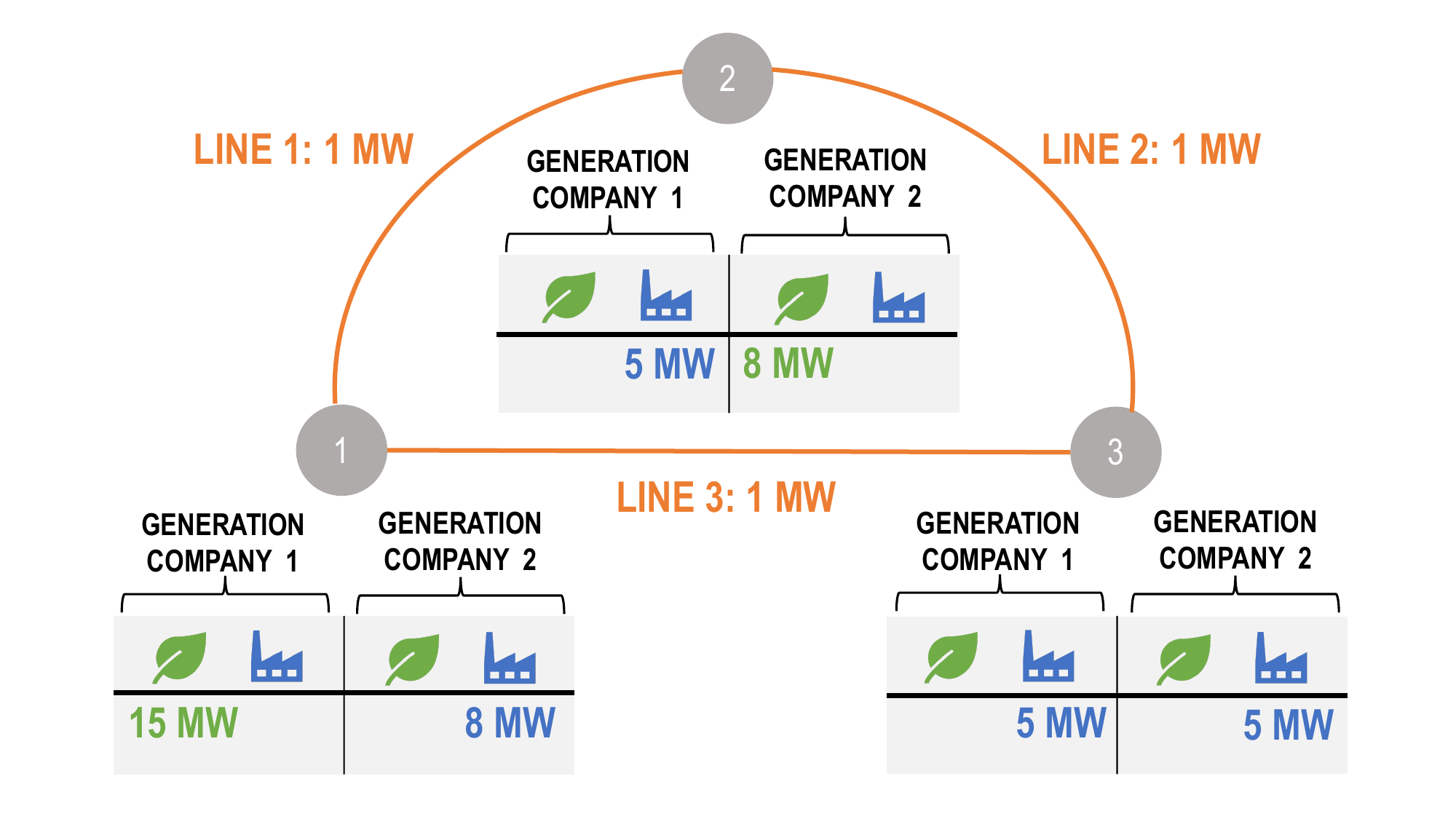}
  \caption{Illustrative energy system structure.}
  \label{fig: illustrative_model_set_up}
 \end{figure}

To be able to conduct a thorough analysis of the optimal investment strategies and generation levels in this section, we first consider a simplified structure for the case study, hereinafter referred to as the illustrative instance. Then in Section \ref{sec: nordic_case}, we consider a case study for the Nordic region. The structure of the illustrative instance is presented in Figure \ref{fig: illustrative_model_set_up}.

The illustrative energy system comprises 3 nodes, 3 transmission lines of a capacity of 1MW each and 2 GenCos that own some conventional and/or renewable power capacities at each of the nodes presented in Figure \ref{fig: illustrative_model_set_up} (green leaf and blue factory images represent VRE and conventional energy sources, respectively). The demand profile at each node, along with the VRE availability factor represented as per cent of the total capacity of the corresponding power plant, are defined in Table \ref{tab: illustrative_example_demand_vres_avail}. 

\begin{table}[h]
\centering
\begin{tabular}{c|c|c|cc}
\multicolumn{1}{c|}{\multirow{2}{*}{Node}} & \multicolumn{1}{c|}{\multirow{2}{*}{VRE availability}} & \multicolumn{1}{c|}{\multirow{2}{*}{Demand profile}} & \multicolumn{2}{c}{Inverse demand function}                 \\ \cline{4-5} 
\multicolumn{1}{c|}{}                      & \multicolumn{1}{c|}{}                                   & \multicolumn{1}{c|}{}                                & \multicolumn{1}{c|}{Slope}  & \multicolumn{1}{c}{Intercept} \\ \hline \hline
1                                          & Moderate (50 \%)                                        & Low                                                  & \multicolumn{1}{l|}{0.04}   & 260                           \\ \hline
2                                          & High (70 \%)                                            & Low                                                  & \multicolumn{1}{l|}{0.04}   & 260                           \\ \hline
3                                          & Low (30 \%)                                             & High                                                 & \multicolumn{1}{l|}{0.0075} & 195                           \\ \hline \hline
\end{tabular}
\caption{Illustrative energy system demand and VRE availability profiles. }
  \label{tab: illustrative_example_demand_vres_avail}
\end{table}

%\begin{figure}[h]
  %\includegraphics[width=\linewidth]{Pictures/illustrative_example_demand_vres_avail.pdf}
 % \caption{Illustrative energy system demand and VREs availability profiles }
  %\label{fig: illustrative_example_demand_vres_avail}
 %\end{figure}

We assume node 3 to have a high-demand profile while the other nodes have low-demand levels. Considering the demand definition through the inverse demand function, the low- or high-demand profiles imply steeper or milder slopes of the function curve, accordingly. The VRE availability factor is moderate and high at nodes 1 and 2, respectively, and low at node 3. For the sake of simplicity, we only assume a single scenario for demand and VRE availability eliminating the stochastic nature of the problem. As for the planning horizon, we consider two consecutive months. To model each month, we multiply hourly data by 720 (i.e., 30 days times 24 hours) to scale it to a month's equivalent. 

Alternative values for the total transmission capacity expansion budget (TEB), GenCos' generation infrastructure expansion budgets (GEB), and carbon tax and renewable generation incentives are considered. The renewable generation incentive is defined as a percentage of total cost reduction arising from installing additional VRE capacity that is subsidised. We performed a univariate sensitivity analysis to understand the impact of each of these mechanisms taken one at a time. Due to a high number of input parameters under analysis, we assume that the input parameters take the discrete values shown in Table \ref{tab: illustrative_case_input_parameters}, where the symbol ``M''(``B'') stands for million (billion). 

\begin{table}[H]
\centering
\begin{tabular}{l|cl}
Input parameter                          & \multicolumn{2}{c}{Discrete values set}                                         \\ \hline \hline
Transmission capacity                    & \multicolumn{2}{c}{\multirow{2}{*}{100K, 1M, 10M, 25M, 50M, 75M, 100M }} \\
expansion budget (\euro)  & \multicolumn{2}{c}{}                                                            \\ \hline
                                         & (GenCo 1: 1M, GenCo 2: 1M),            & \multirow{2}{*}{\Big\}small}            \\
Generation infrastructure                & (GenCo 1: 10M, GenCo 2: 10M)             &                                     \\ \cline{2-3} 
expansion budget (\euro)  & (GenCo 1: 10M, GenCo 2: 1B),              & \multirow{2}{*}{\Big\}large}            \\
                                         & (GenCo 1: 1B, GenCo 2: 1B)              &                                     \\ \hline
Carbon tax (\euro  / MWh) & \multicolumn{2}{c}{0, 25, 50, 75, 100, 125, 150, 175, 200}                  \\ \hline
Renewable generation incentives (\%)     & \multicolumn{2}{c}{0, 10, 20, 30, 40, 50, 60, 70, 80, 90, 100}              \\ \hline \hline
\end{tabular}
\caption{Illustrative case study input parameters.} 
\label{tab: illustrative_case_input_parameters}
\end{table}

For each combination of the input parameters, we analyse the total welfare (the value of the upper-level problem objective function), the share of VRE in the total generation mix and the total generation level, which we refer to as output factors. %The primary purpose of the last output factor (total generation level) is to identify whether changing the input parameters values impacts the change in the total generation level. 

\subsection{Summary of the results}
The results for this illustrative case study have proven ineffective to simply increase the TEB in case the GenCos do not plan significant expenses to expand generation infrastructure, i.e., more than \euro1B. However, in case such circumstances occur, increasing the TEB leads to an increase in the share of VRE in the total generation mix without dampening the total amount of energy generated. However, the question regarding the effectiveness of such a policy in terms of the ratio between its cost and impact remains. 

When VRE generation capacity expansion incentive is the single VRE share-increasing policy being considered, one can notice that the numerical results suggest unstable behaviour in the values of each of the output factors. Depending on the GenCos GEB, the share of the subsidised costs might be required to be significantly high to demonstrate a positive correlation between the VRE share and total welfare. Moreover, while the VRE share and total welfare would experience an increase, the total generation amount might decrease before it starts rising again when the incentive increases (see Figures \ref{fig: illustrative_case_subsidy_central} and \ref{fig: illustrative_case_subsidy_perfect}). Hence, the policymaker should carefully consider at which level the incentives will not harm total generation. Furthermore, it is essential that policymakers consider the increasing global demand for electricity attributed to the digitalisation and electrification of diverse domains, alongside the expanding industrial output and services sector \cite{mir2020review}.

Applying a carbon tax as a unique renewables-driven policy expectedly allows for an increase in VRE share, which induces a decline in the amount of energy produced and a decrease in total welfare value. However, depending on GenCos' GEB, there might be a threshold in the carbon tax value only after which the anticipated rise of VRE share will occur. Nevertheless, it is worth mentioning that even in the absence of a carbon tax or VRE generation capacity expansion incentives there is a consistent trend favouring renewable energy technologies regarding levelised costs of electricity generation, as compared to power production based on fossil fuels \cite{timilsina2022economics}.

It is worth highlighting that conducting sensitivity analysis for only one parameter at a time allowed us to obtain detailed information regarding all the phenomena that occurred and significantly simplified the interpretation of the results. Introducing multi-dimensional parameters in sensitivity analysis drastically reduces the possibilities for meaningful inferences. Nevertheless, it allows for investigating the influence of the composition of the input parameters on the output factors. This analysis is conducted in the next Section \ref{sec: nordic_case} while applying the proposed method to a realistic case study.

The above results come from the analysis presented in Sections \ref{sec: il_TEB}, \ref{sec: il_VRE} and \ref{sec: il_tax}. It is worth highlighting that in what follows the results for the markets with the centralised planner and perfect competition are nearly identical. Such phenomenon can be also observed in \citep{VIRASJOKI2020104716}.

\subsection{Transmission infrastructure expansion budget (TEB sensitivity)}
\label{sec: il_TEB}
The proposed methodology is applied to derive an optimal infrastructure expansion strategy and generation levels considering carbon tax and renewable subsidies to be fixed to 0 \euro \ / MWh and 0\%, respectively, while the GEB takes discrete values given in Table \ref{tab: illustrative_case_input_parameters}. Meanwhile, we conduct a sensitivity analysis on the TEB considering values given in Table \ref{tab: illustrative_case_input_parameters}. The experiments are replicated for both a centralised planner and perfectly competitive generation markets (\eqref{slm: objective}--\eqref{slm: flow_primal_feas} and \eqref{bi_level_formulation}, respectively). 

\begin{figure}[h!]
  \includegraphics[width=\linewidth]{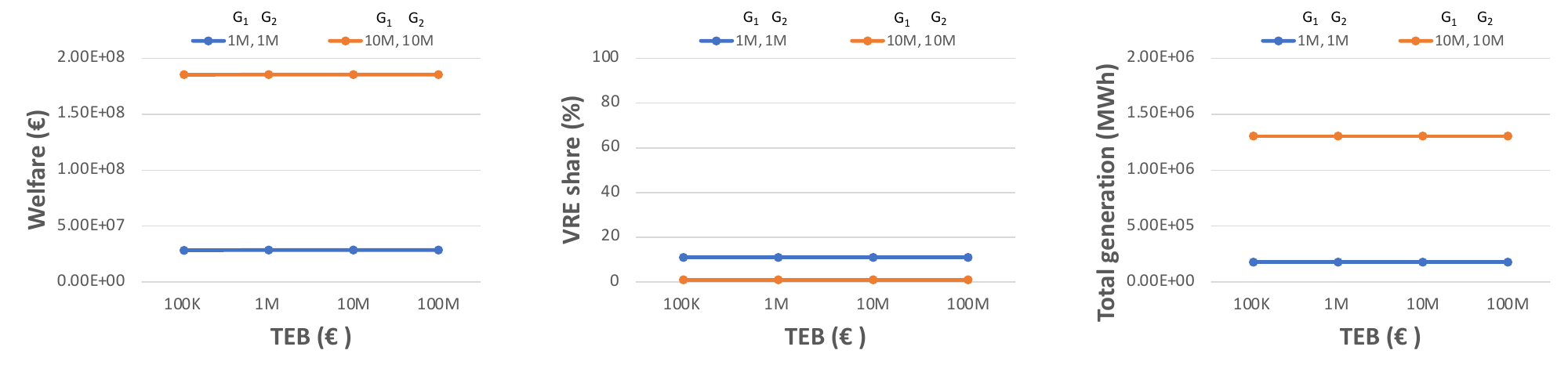}
  \caption{Sensitivity analysis on the TEB considering the small GEB under centralised and perfect competition settings.}
  \label{fig: illustrative_case_transm_small_budget}
\end{figure}

Figure \ref{fig: illustrative_case_transm_small_budget} illustrates the results of the sensitivity analysis on the TEB considering the lowest levels of GEB, i.e., 1M \euro \ or 10M \euro \ for each of the GenCos, indicated by blue and orange lines, respectively. In the graph, each line's legend contains two values corresponding to the GEB values of GenCo 1 (G$_1$) and GenCo 2 (G$_2$), respectively. The results suggest no influence of the TEB value on any of the output factors. Moreover, the optimal investment and generation levels are identical for both centralised and perfectly competitive market structures, which is due to the GEB-associated constraint being binding. 

\begin{figure}[h!]
  \includegraphics[width=\linewidth]{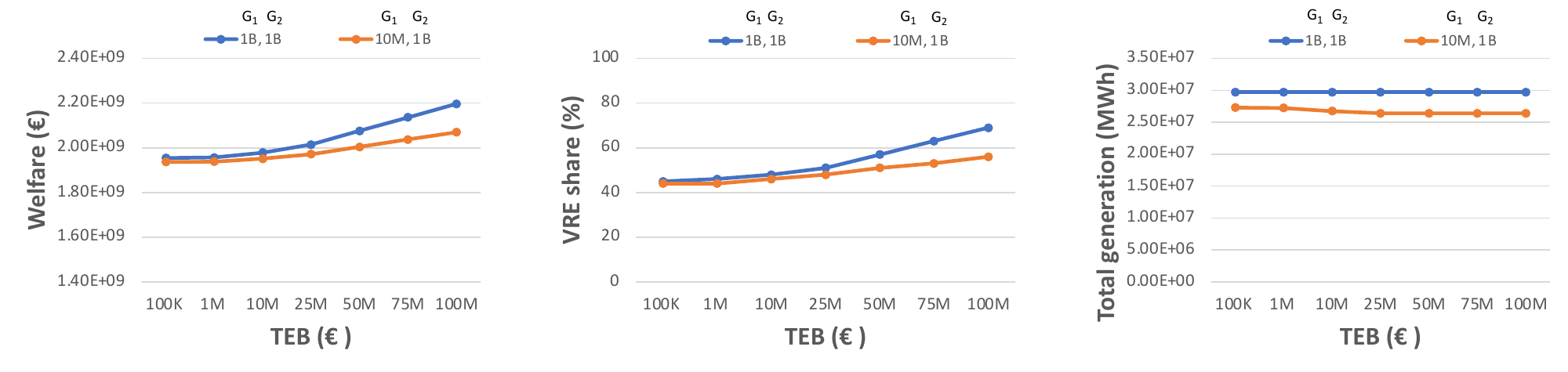}
  \caption{Sensitivity analysis on the TEB considering the large GEB under centralised planner.}
  \label{fig: illustrative_case_transm_large_budget_central}
\end{figure}

\begin{figure}[h!]
  \includegraphics[width=\linewidth]{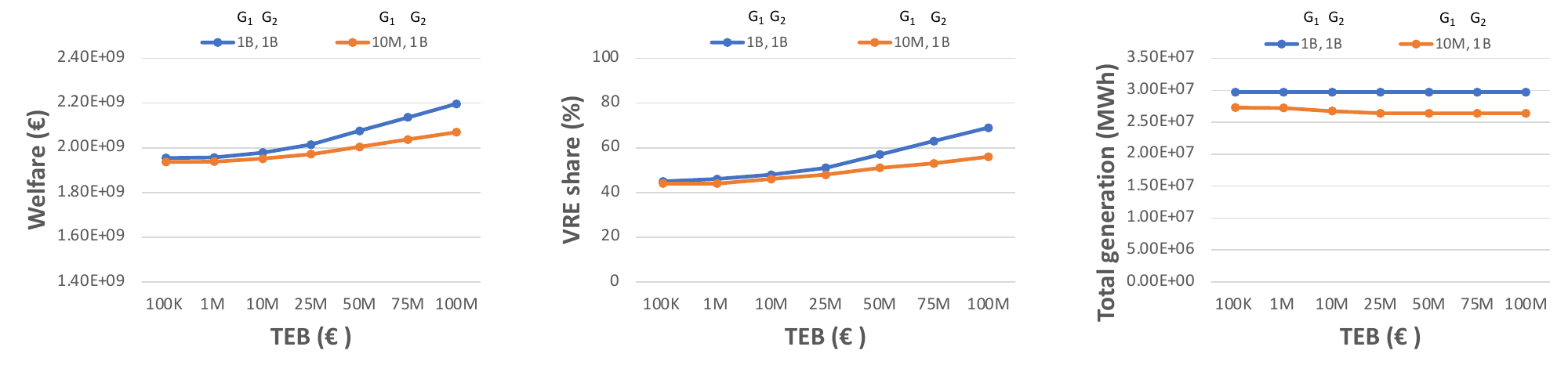}
  \caption{Sensitivity analysis on the TEB considering the large GEB under the perfectly competitive market setting.}
  \label{fig: illustrative_case_transm_large_budget_perfect}
\end{figure}

Figures \ref{fig: illustrative_case_transm_large_budget_central}  and \ref{fig: illustrative_case_transm_large_budget_perfect} illustrate the influence of changing the TEB value on each of the output factors assuming large GEB values and considering centralised and perfect competition settings, respectively. 

 \begin{figure}[h!]
  \includegraphics[width=0.8\linewidth]{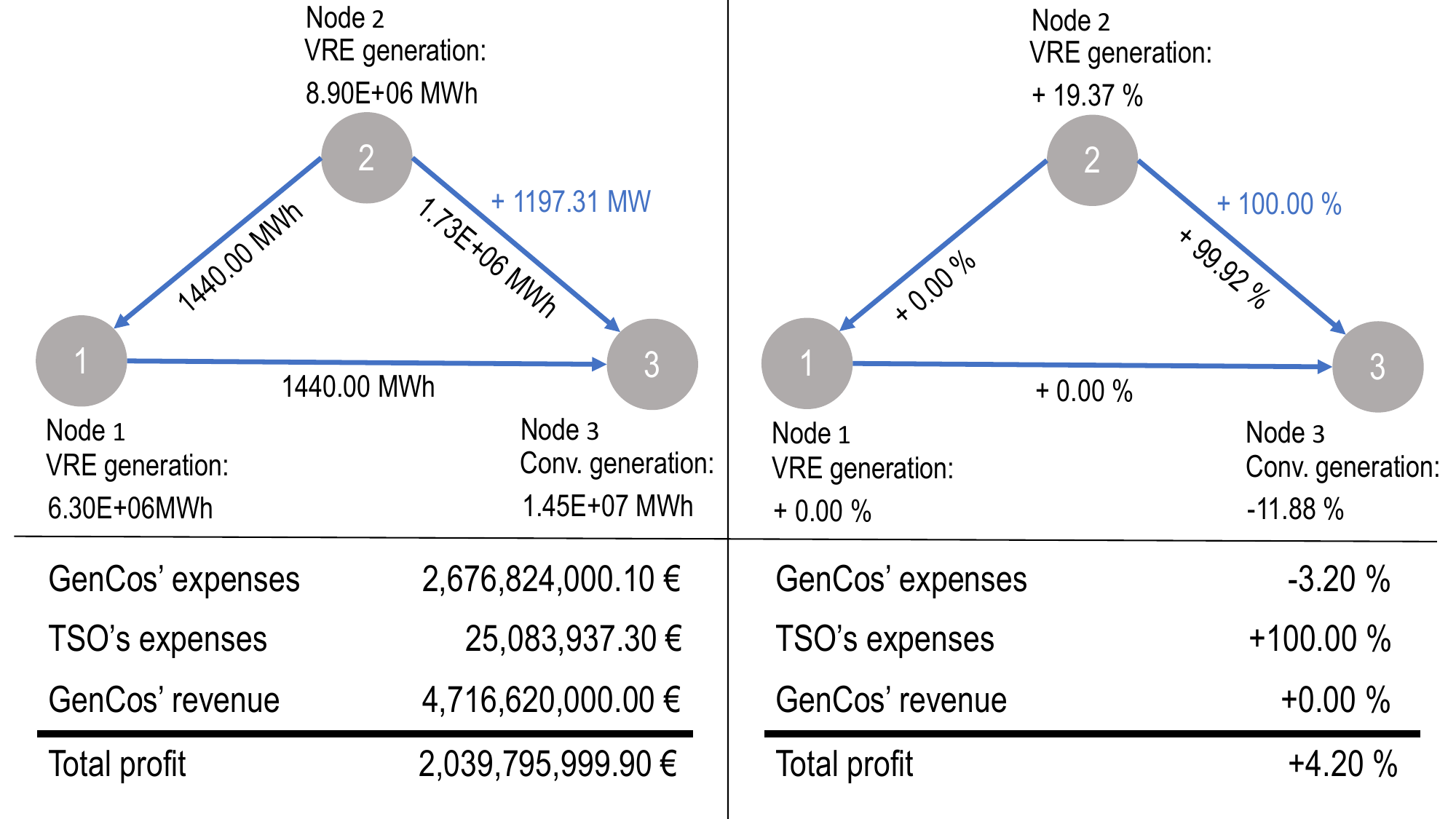}
  \centering
  \caption{Differences when considering no incentives, no carbon tax and \euro1B GEB for each GenCo but changing TEB from \euro25M (left frame ) to \euro50M (right frame) in perfectly competitive market.}
  \label{fig: illustrative_case_transm_diff_1}
\end{figure}

 As can be seen, larger TEB values lead to an increase in the renewable generation share along with the total welfare value, and the effect is magnified with the increase of the GEB. Apparently, there is a very small impact on the total generation. Hence, the sole increase of the transmission capacity without involving other policies leads to substitution for lower-priced VRE generation followed by an increase of the total welfare with a slight decline in the total generation.

Let us consider the case of a \euro1B GEB for each GenCo and compare the differences in optimal investment and generation portfolios when increasing TEB from \euro25M to \euro50M in a perfectly competitive market. The detailed analysis for this particular example is presented in Figure \ref{fig: illustrative_case_transm_diff_1}. The arrow in the figure indicates the direction of the power flow. The number with the sign ``+''  next to it on the left represents the capacity added to the transmission line (in MW). And the number without a sign indicates the total amount of energy transmitted during all time periods through this line (in MWh). The right frame of Figure \ref{fig: illustrative_case_transm_diff_1} demonstrates the relative change in the output factors' values in relation to the values presented on the left. As one can notice, doubling the TEB from \euro25M to \euro50M is followed by an increase of 100\% in capacity and 99.92\% in the flow through the line connecting node 2 (with the highest VRE availability) and node 3 (with the highest demand profile). This transmission capacity expansion motivates GenCos to invest in more VRE generation at node 2, increasing the VRE generation by 19.37\% while reducing the conventional generation at node 3 by 11.88\%. Ultimately, this new configuration leads to the decrease of GenCos' costs by 3.20\% without any decrease in total generation levels. The latter phenomenon, in turn, leads to an increase in the profit by 4.20\% and, hence, the increase in the total welfare as well. This illustrates that, even in this stylised example, the model is capable of capturing the key features of the problem regarding the availability of system connectivity and its effect on Gencos' motivation to expand their VRE generation.

\subsection{VRE generation incentives}
\label{sec: il_VRE}
Next, the VRE generation-capacity subsidies are varied. Similarly to the previous subsection, we consider the values of the remaining parameters to be fixed. In particular, we assume the carbon tax to remain at 0 \euro \ / MWh, the TEB at \euro10M while considering the GEB values as defined in Table \ref{tab: illustrative_case_input_parameters} but excluding the case where both GenCos have \euro1M GEB, due to the similarities in optimal decisions values with the case in which they have the same budget set to \euro10M.

\begin{figure}[h!]
  \includegraphics[width=\linewidth]{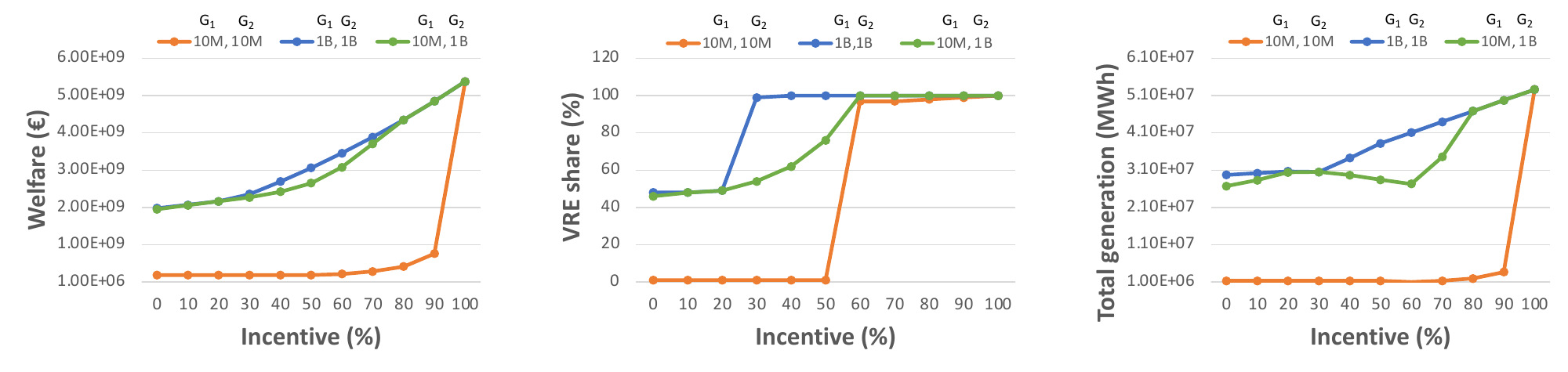}
  \caption{Sensitivity analysis on the VRE generation capacity subsidies considering different GEB values under centralised planner.}
  \label{fig: illustrative_case_subsidy_central}
\end{figure}

\begin{figure}[h!]
  \centering
  \includegraphics[width=\linewidth]{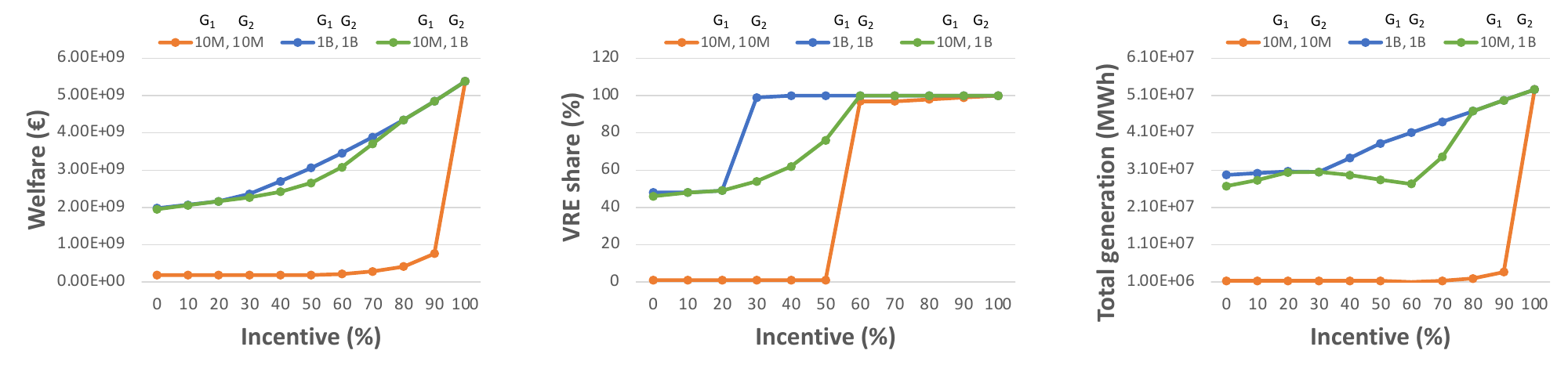}
  \caption{Sensitivity analysis on the VRE subsidies considering different GEB values in perfect competition}
  \label{fig: illustrative_case_subsidy_perfect}
\end{figure}

Figures \ref{fig: illustrative_case_subsidy_central} and \ref{fig: illustrative_case_subsidy_perfect} summarise the values of the output factors considering different levels of VRE generation capacity subsidies for the centralised planner and perfectly competitive market, respectively.

The numerical experiments suggest that in cases where the GenCos possess identical GEBs, there is a threshold in the percentage value of VRE subsidies after which %the marginal cost of producing additional units of VRE energy becomes less expensive. The latter makes 
the decision to invest in VRE energy becomes more appealing (i.e., more profitable) for the GenCos, when compared to conventional generation sources. In particular, for this illustrative example, this value lies between 50\% and 60\% in the case of lower GEB, i.e., \euro10M for each GenCo, and between 20\% and 30\% when this budget is increased to \euro1B for each GenCo. 

\begin{figure}[h!]
  \centering
  \includegraphics[width=0.8\linewidth]{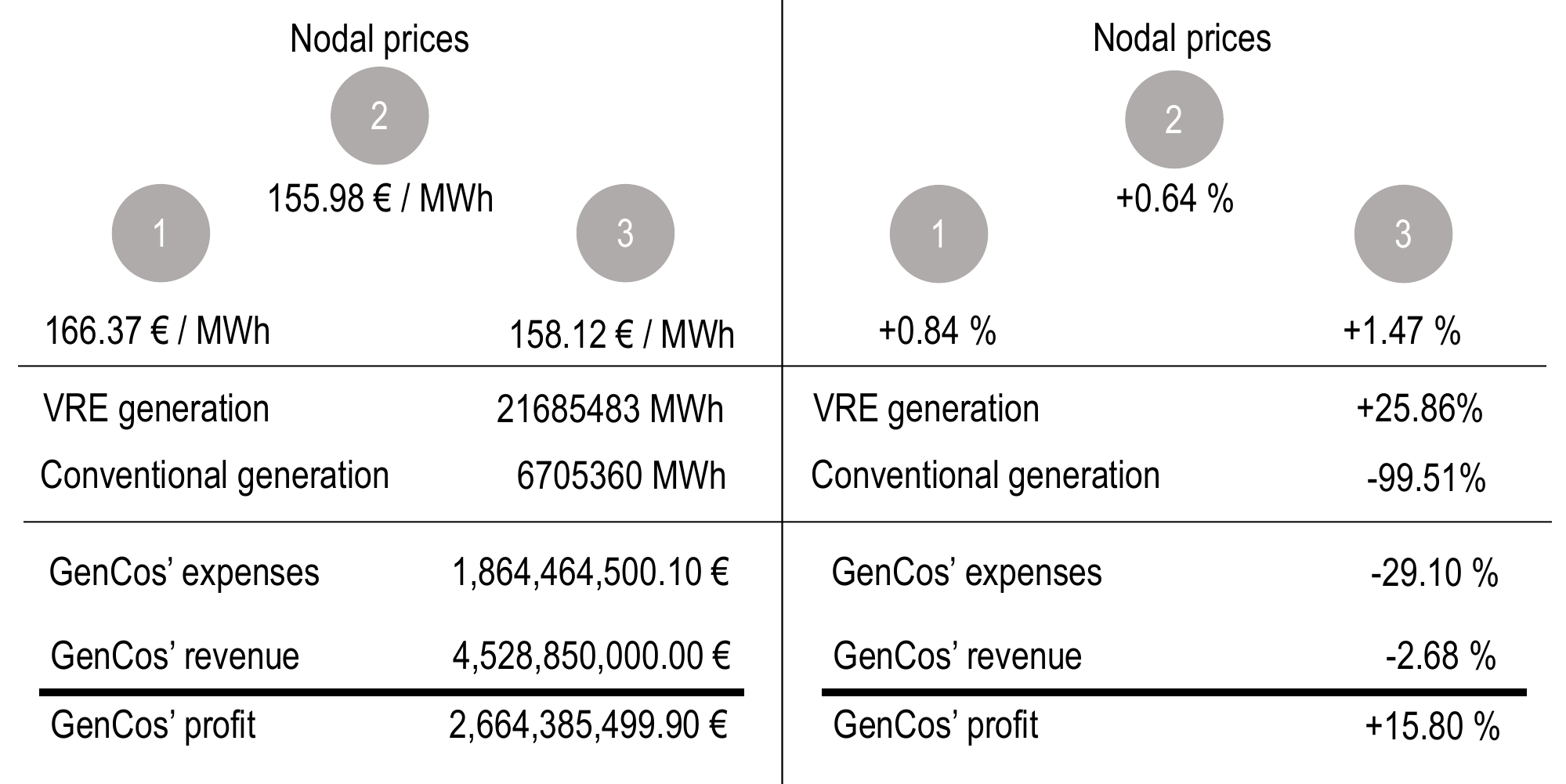}
  \caption{Differences in optimal results when changing VRE incentive from 50\% (left-hand side) to 60\% (right-hand side).}
  \label{fig: illustrative_case_incentive_central_diff}
\end{figure}

Another unanticipated phenomenon observed is the loss in the total generation that occurs when subsidies are between 50\% and 60\% (30\% and 60\%; 20\% and 30\%) of the VRE generation expansion costs for the GenCos with \euro10M both (\euro10M and \euro1B; \euro1B both) GEB. This is caused by a combination of higher investments in VRE generation capacity and lower generation levels when compared to the case without subsidies. This, in turn, is due to two phenomena. Firstly, the inverse linear dependence between demand and price allows for higher prices when generation levels are lower. Hence, the loss in GenCos' revenues at each node (nodal revenues) is minimised. Secondly, the GenCos experience significant savings in expenses, which compensate the decrease in nodal revenues. For example, Figure \ref{fig: illustrative_case_incentive_central_diff} presents a detailed comparison between optimal investment decisions and generation levels considering 50\% and 60\% incentives for GenCo 1 and GenCo 2 with \euro10M and \euro1B GEB, respectively, in case of perfect competition. As one can notice, the increase of subsidies, in this case, stimulates a decrease in conventional generation by almost 100\%, which is not completely substituted by VRE generation. However, the fall in GenCos' revenue is compensated by increased nodal prices and overcompensated through reductions in GenCos' expenses by 29.10\%, ultimately leading to a profit increase of 15.80\%. 

\subsection{Carbon tax}
\label{sec: il_tax}
Lastly, we present a sensitivity analysis for the carbon tax. Similarly to the case of the other input parameters, we fix the values of incentives and TEB budget to 0\% and \euro10M, respectively, while assuming the GEB to take values from Table \ref{tab: illustrative_case_input_parameters}. Here, we have also omitted the cases in which GenCos have a \euro 1M capacity expansion budget due to the similarities in optimal decision values with the \euro10M case. We consider the values for the carbon tax to lie within the interval given in Table \ref{tab: illustrative_case_input_parameters}. 

\begin{figure}[h!]
  \includegraphics[width=\linewidth]{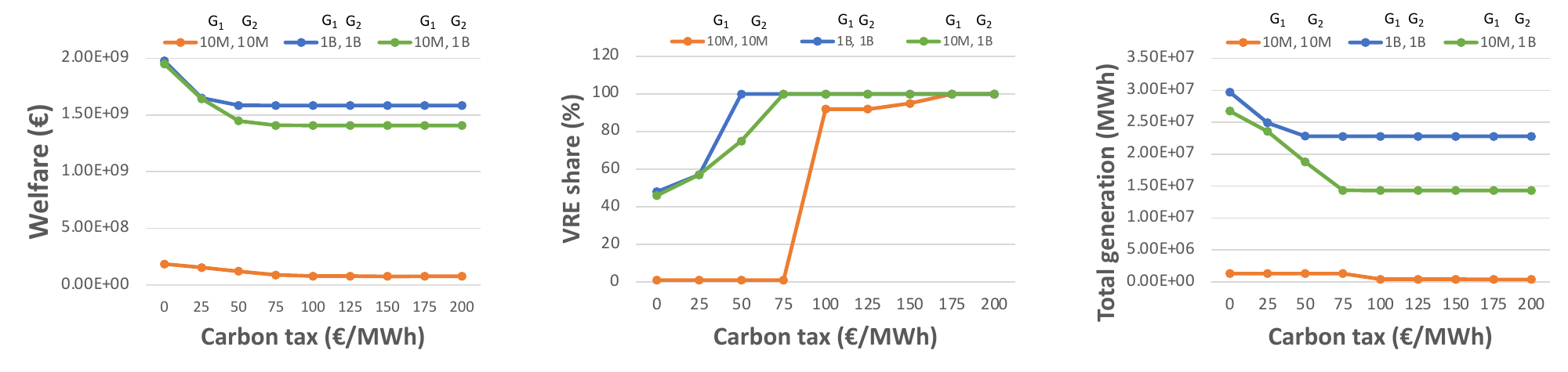}
  \caption{Sensitivity analysis on the carbon tax considering different budget values for generation capacity expansion under centralised planner}
  \label{fig: illustrative_case_tax_central}
\end{figure}

\begin{figure}[h!]
  \includegraphics[width=\linewidth]{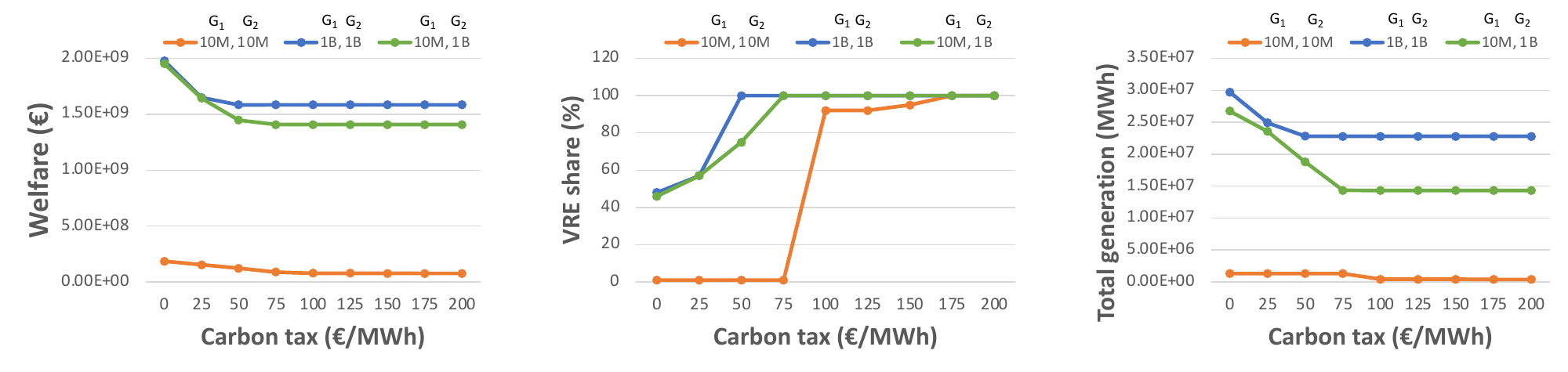}
  \caption{Sensitivity analysis on the carbon tax considering different GEB values under perfect competition}
  \label{fig: illustrative_case_tax_perfect}
\end{figure}

Figures \ref{fig: illustrative_case_tax_central} and \ref{fig: illustrative_case_tax_perfect} present the output factors considering centralised and perfectly competitive market structures, respectively. As one can see, the results present a threshold value for the tax, after which it becomes effective and influences the increase of the VRE share in the total mix. This value lies between 75 \euro \ / MWh and 100 \euro \ / MWh for the GenCos with \euro10M GEB, and between 0 \euro \ / MWh and 25 \euro \ / MWh when one or both GenCos posses \euro1B GEB. Nevertheless, while serving its purpose regardless of the GEB value, the carbon tax causes a decrease in the total generation and, consequently, in the total welfare. Both output factors only remain stable once GenCos cease nearly all the conventional generation (i.e., VRE share in the total generation mix is close to 100\%). 

\section{Numerical experiments: The Nordic energy system case study}
\label{sec: nordic_case}
\begin{figure}[h]
\centering
  \includegraphics[width=0.50\linewidth]{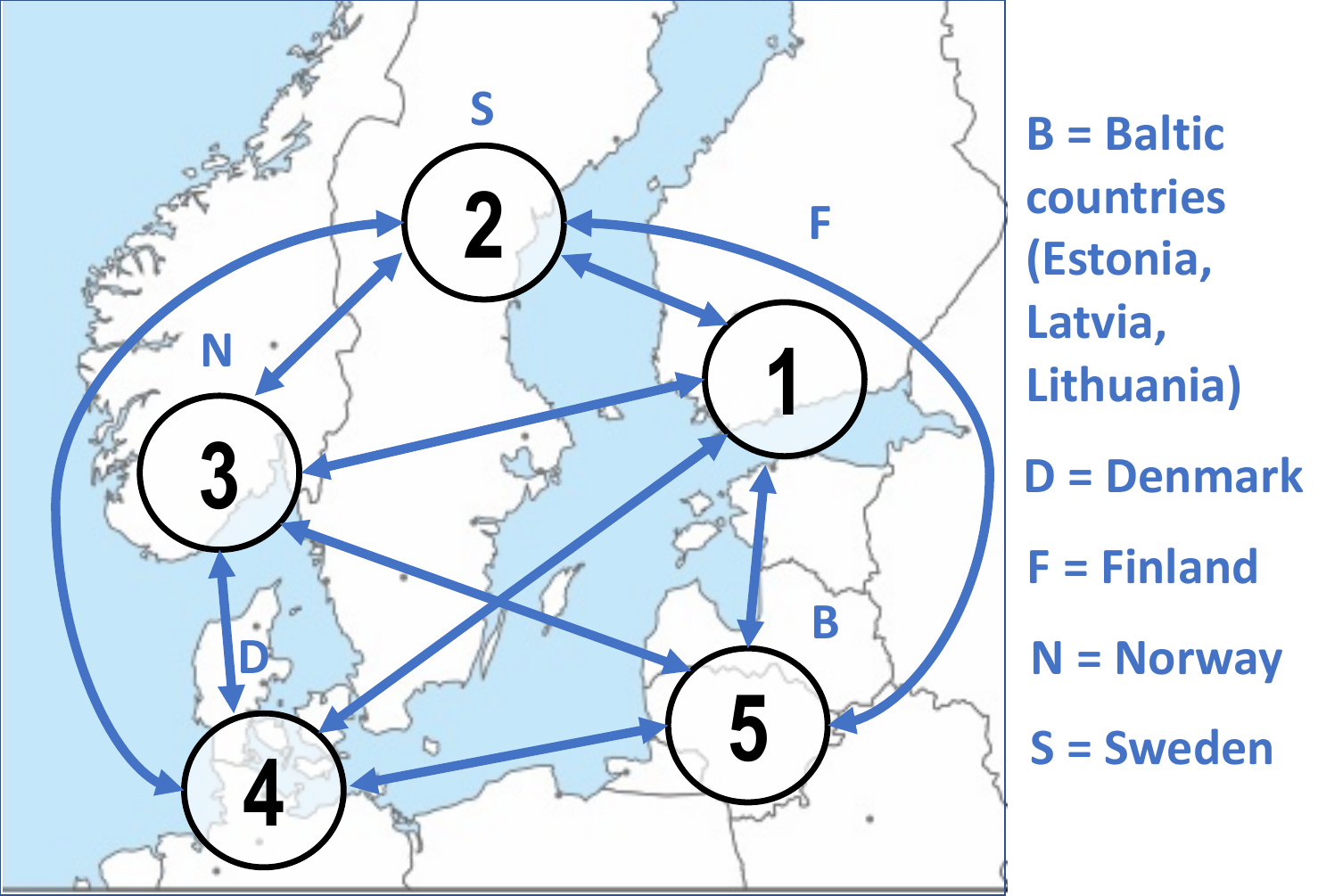}
  \caption{Case study: Nordics' energy system \citep{noauthor_scandinavia_nodate}.}
  \label{fig: nordics_map}
 \end{figure}

Next, we apply the proposed modelling assessment to a system representing the features of the Nordic energy system, enlarged with the Baltic countries grouped as one node. Therefore, the case study system involves 5 nodes representing Finland, Sweden, Norway and Denmark (Nordic countries) and consolidated Estonia, Latvia and Lithuania (Baltic countries), assuming all possible interconnections between these nodes to be potential or existing transmission lines. We conduct an analysis analogous to Section \ref{sec: illustrative_example} but simultaneously changing different input parameters. The schematic structure of the case study energy system is illustrated in Figure \ref{fig: nordics_map}.

For the purpose of simplification, we consider the existence of 5 GenCos, where each is associated with a distinct node (as shown in Figure \ref{fig: nordics_map}) and possesses all the generation capacity installed at the corresponding node prior to the beginning of the modelling horizon. This implies that none of the GenCos is assumed to own any preinstalled generation capacity at other nodes. Nevertheless, we assume that each of the GenCos may invest in the generation capacity at any node. This assumption is based on the illustrative fact that Finnish energy generation company Fortum operates in several countries, such as Finland, Sweden, Norway, Poland, Germany, Denmark, India and others \citep{fortum_country_nodate}. Additionally, we enlarge the set of conventional energy sources by considering nuclear energy, open-cycle and combined-cycle gas turbines (OCGT and CCGT, respectively), coal and biomass. As VRE sources, we consider solar as well as onshore and offshore wind. We also account for the hydropower capacity in these countries, though we assume they can not be expanded.

\subsection{Data}
All the data for demand, energy prices, installed generation capacity, VRE availability and transmission parameters were gathered from various sources, with information ranging from the years 2018 to 2020. One can refer to Appendix  \ref{appendix: data} for more information regarding the data.

\paragraph{Input parameters}
Due to the large-scale nature of the case study instance and the consequent computational challenges we have faced, we limited the range of the values for each of the input parameters to a discrete set with two values in the sensitivity analysis. The first one represents a possible ``low'' value, and the second one represents a ``high'' value for a corresponding parameter. Therefore, we solve the case study instance for each possible set of 4 input parameters where each parameter takes one of two possible values from a corresponding set which leads to solving 16 different instances in total. The values for the parameters are presented in Table \ref{Tab: Nordic_case_input_parameters}. It is worth mentioning that the data for the TEB and GEB is presented for an annual time frame. More information on the procedure for generation input parameter values can be found in Appendix \ref{appendix: data}.

\begin{table}[H]
\centering
\begin{tabular}{l|c|c}
Input parameter & ``Low'' value  & ``High'' value \\ \hline  \hline 
{Transmission capacity} & \multirow{2}{*}{625M } &  \multirow{2}{*}{1.25B}    \\ 
{expansion budget (\euro / year)} & &    \\ \hline
{Generation capacity} & \multirow{2}{*}{3B} &  \multirow{2}{*}{6B }    \\ 
{expansion budget (per GenCo, \euro / year)} & &    \\ \hline
 Carbon taxes (\euro  / ton of CO$_2$) & 8 & 71    \\ \hline
Renewable generation incentives (\%) & 5 & 20
\\ \hline \hline  
\end{tabular}
\caption{Discrete values for the input parameters for the large-scale case study.}
\label{Tab: Nordic_case_input_parameters}
\end{table}

\paragraph{Numerical issues with large-scale MPEC}
Despite being able to solve the illustrative instance presented in Section \ref{sec: illustrative_example}, the state-of-the-art mathematical optimization solvers we tested were unable to solve this more realistic instance. We experimented with different complementary slackness reformulation options and primal-dual equality reformulation available in the BilevelJuMP.jl package \citep{Joaquim_bilevel_2022} along with reductions in the instance size. However, trying to solve a simplified instance with Gurobi v.9.5.0 \citep{gurobi} or CPLEX v.12.10 \citep{noauthor_ibm_2021} within 48 hours did not lead to a solution with integrality gaps lower than 100 \%. An analogous issue has been faced by \citet{VIRASJOKI2020104716} where the authors were not able to solve their large-scale mixed-integer quadratically constrained quadratic programming problems with any of the state-of-the-art solvers available in GAMS \citep{noauthor_gams_nodate}. 

Therefore, following the solution approach proposed in \citep{VIRASJOKI2020104716}, we applied an iterative procedure in which we discretise and exhaustively enumerate and fix all possible upper-level decisions and solve the lower-level problem. Then, we determine ex post the decisions that yield optimal solutions. The discrete values we considered for the capacity expansion of the transmission lines are \{0 MW, 3000 MW, 6000 MW, 9000 MW\}. Taking into account the existence of 10 transmission lines in the system, the aforementioned enumerative set leads to 1,048,576 lower-level problems that must be solved for each combination of the input parameters. 

\paragraph{Representative days}
In light of the large number of lower-level problems to be solved for each set of input parameter values, we considered 3 scenarios and 4 time periods. Hence, each scenario spans a group of 3 representative days corresponding to the demand profile, solar and wind availability, respectively. And the day (24 hours) timeline is equally divided into four 6-hour time periods. The idea is inspired by \citet{LI2018965}, who used 16 time slots (four diurnal times-slices in four seasons) to model a year for the energy system. \citet{LIND201744} used a similar approach for TIMES-Oslo and TIMES-NORWAY models. Additionally, \citet{FODSTAD2022112246} mentioned other authors who used representative daily time slots as a part of the modelling. 

To derive representative days, we used hierarchical clustering \citep{Teichgraeber2019joss} for the hourly data (8760 h) spanning the time series for demand, wind and solar PV power production for each of the nodes. We highlight that all the data has been normalised, implying that each element of the 3 data sets obtained a value between 0 and 1, where 1 refers to the highest regional demand in the case of the demand-related data. For the solar PV and wind production data sets, the normalisation values are based on the installed capacities, and hence the assigned normalised values never reach 1. When conducting the hierarchical clustering, no weights have been considered between the data time series. The representative days for each cluster have been selected for each data set to most closely approximate the medoid of a cluster. The allocation of weights between the clusters, as well as the distribution of representative days for each cluster and node, are presented in Appendix \ref{ap: hierarchical clustering}. For example, scenario one for Finland has a weight of $\frac{145}{365}$ and is represented by a set of days $\{43, 12, 12\}$ corresponding demand profile, wind and solar availability, respectively. It is important to highlight that the notions of clusters and the weights associated with clusters are used in this context to represent the scenarios and the scenario probabilities, respectively, denoted by $s \in S$ and $P^s, \forall s \in S$. 

\subsection{Numerical results}
\begin{figure}[h!]
\centering
  \includegraphics[width=\linewidth]{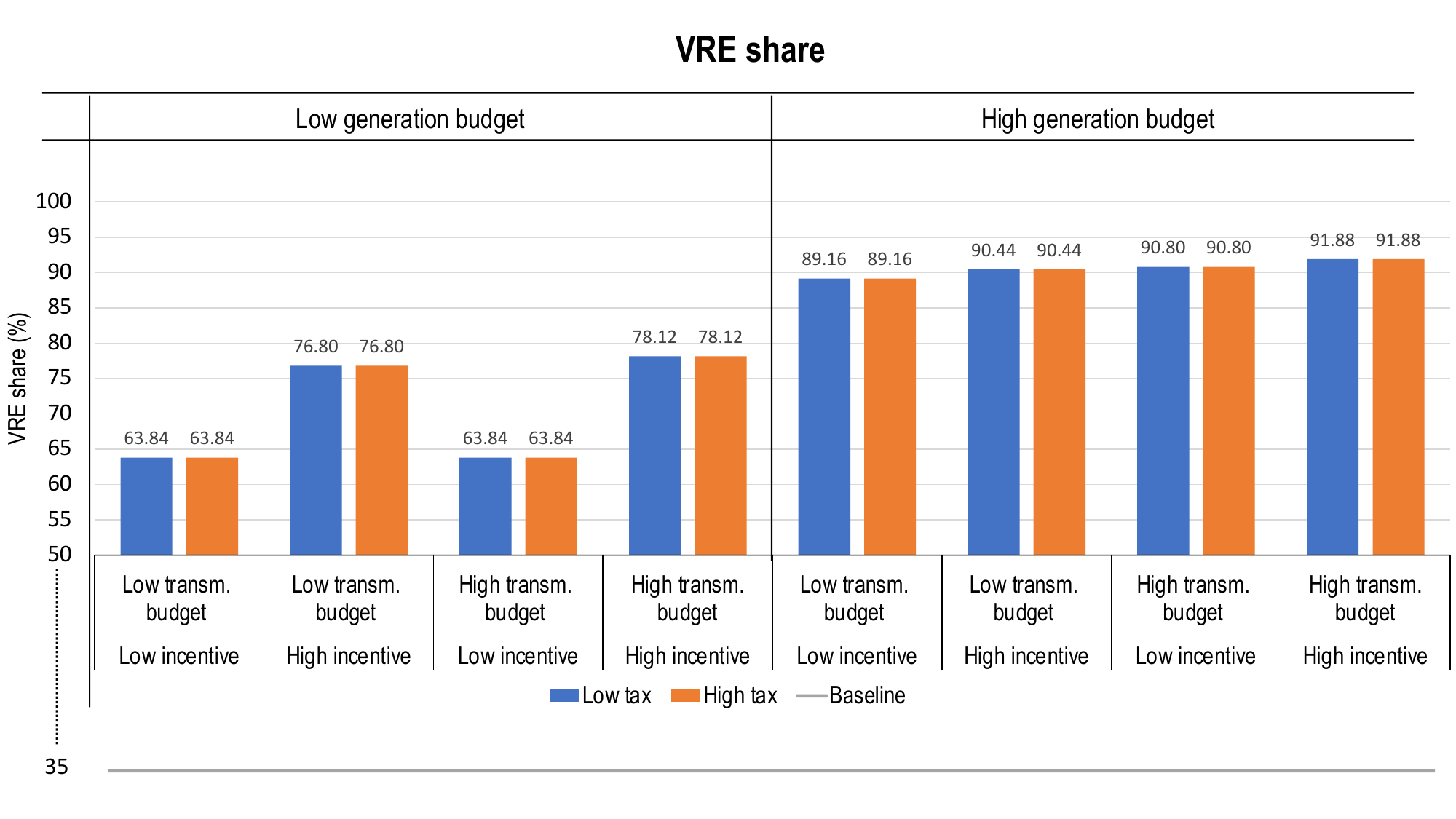}
  \caption{Case study: Impact of carbon taxes on VRE share.}
  \label{fig: nordics_tax_vres}
 \end{figure}

\begin{figure}[h!]
\centering
  \includegraphics[width=\linewidth]{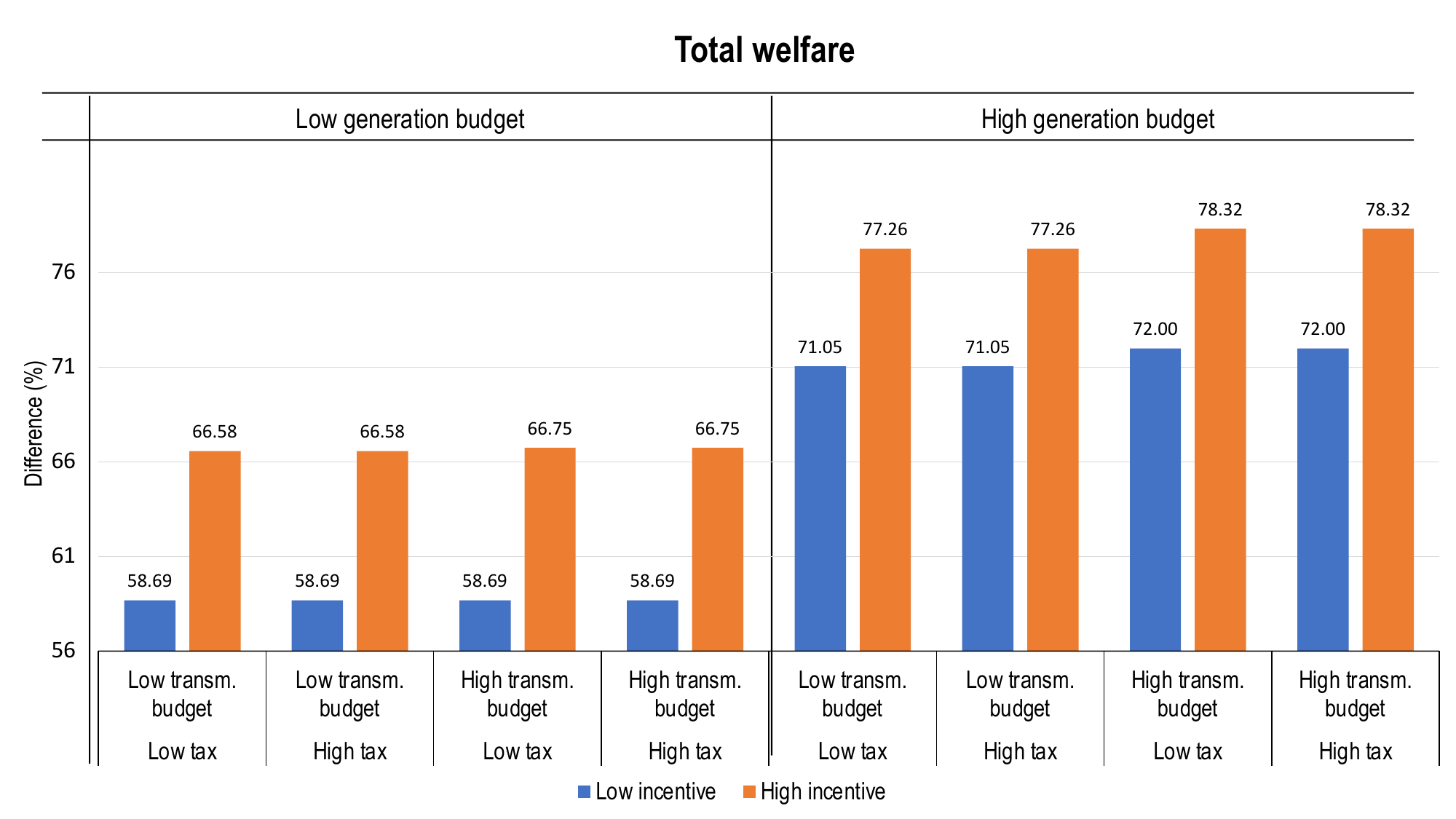}
    \caption{Case study: Impact of VRE investment incentive on welfare.}
  \label{fig: nordics_incentive_welfare}
 \end{figure}

\begin{figure}[h!]
\centering
  \includegraphics[width=\linewidth]{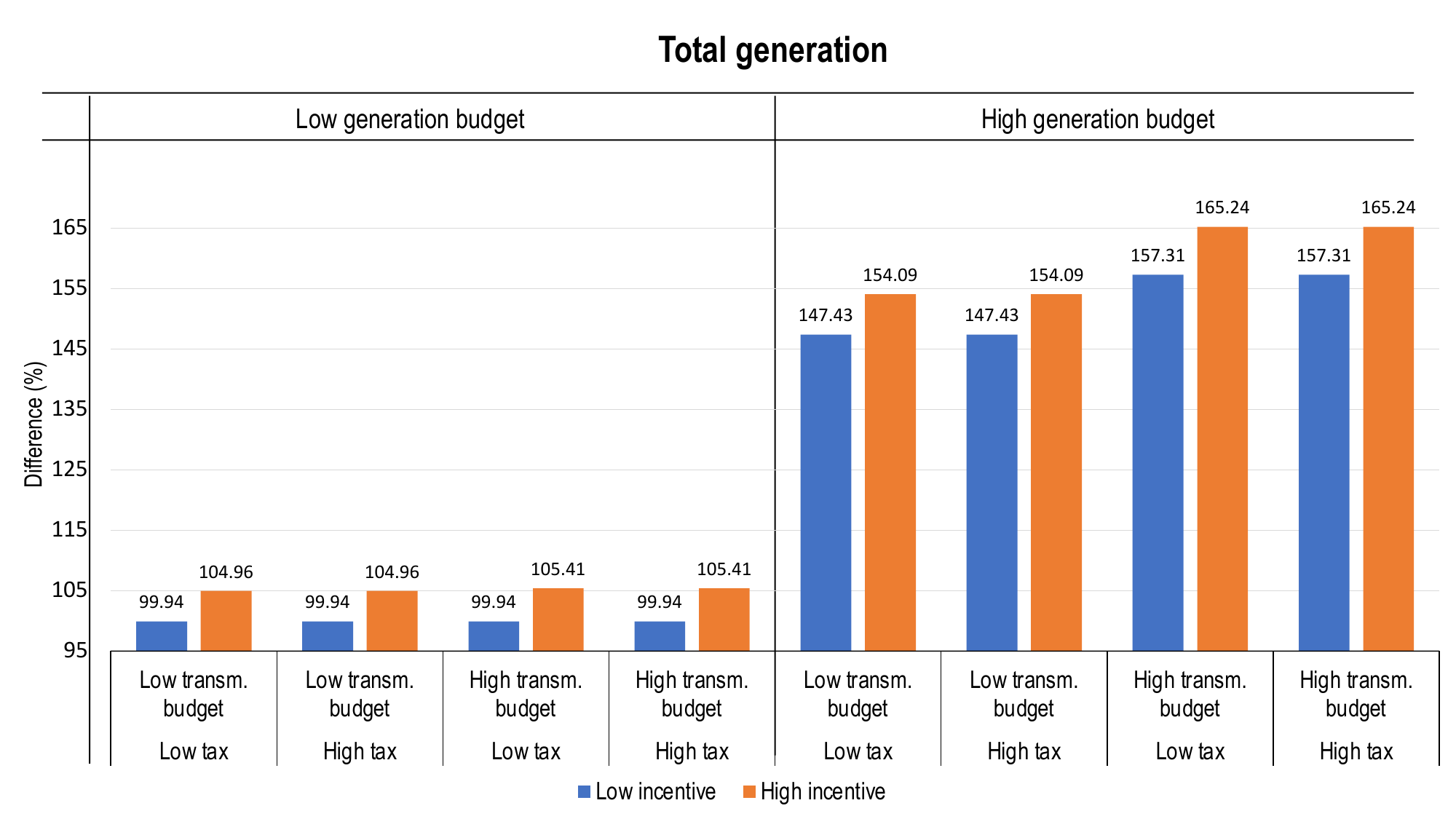}
   \caption{Case study: Impact of VRE investment incentive on total generation.}
  \label{fig: nordics_incentive_generation}
 \end{figure}

Table \ref{tab: nordic_summary} presents the summary regarding the compositions of the input parameter values that have the strongest and the weakest effect on each of the output factors (VRE share in the generation mix, total welfare and generation amount). Additionally, it provides the relative values indicating the improvement of the corresponding factor compared to its optimal value in case all the input parameters, i.e., carbon tax, VRE incentives, TEB, and GEB are zero (baseline case). It is worth highlighting that the ``baseline'' case does not closely reflect reality, as these input parameters are greater than zero in practice. Figure \ref{fig: nordics_tax_vres} presents the optimal share of VRE in the generation mix for different combinations of input factors. The bars indicate the share of VRE when the parameter values are ``low'' or ``high'' according to Table \ref{Tab: Nordic_case_input_parameters}. The ``baseline'' shows the value of VRE in the optimal generation mix when maximising the total welfare and considering all the input parameters to be zero. Figures \ref{fig: nordics_incentive_welfare} and \ref{fig: nordics_incentive_generation} present relative differences in total welfare and generation values, respectively, for different compositions of the input parameters when compared to the case when all of the parameters are set to zero, i.e., no policies are involved. It is worth highlighting that while all the information regarding output results can be obtained from Figures \ref{fig: nordics_tax_vres}, \ref{fig: nordics_incentive_welfare} and \ref{fig: nordics_incentive_generation}, the primary purpose of these figures is to demonstrate the change in output factors values considering different values of distinct input parameters. In particular, in the case of Figure \ref{fig: nordics_tax_vres}, this parameter is a tax value, while Figures \ref{fig: nordics_incentive_welfare} and \ref{fig: nordics_incentive_generation} primarily demonstrate information regarding incentive values. The remaining alternative visualisations of the output data is presented in Appendix \ref{ap: alternative_representations}.

\begin{table}[h]

\begin{tabular}{lllllll}
\hline \hline
\multicolumn{7}{c}{Small GEB}                                                                                                                                                     \\ \hline
\multicolumn{1}{l|}{}                & \multicolumn{2}{c|}{VRES}                        & \multicolumn{2}{c|}{Welfare}                    & \multicolumn{2}{c}{Generation amount} \\ \hline
\multicolumn{1}{c|}{}                & Low tax        & \multicolumn{1}{l|}{}           & Low tax        & \multicolumn{1}{l|}{}          & Low tax             &                 \\
\multicolumn{1}{l|}{Largest effect}  & High incentive & \multicolumn{1}{l|}{+ 42.28\%} & High incentive & \multicolumn{1}{l|}{+ 66.75 \%} & High incentive      & + 105.41 \%     \\
\multicolumn{1}{l|}{}                & High TEB       & \multicolumn{1}{l|}{}           & High TEB       & \multicolumn{1}{l|}{}          & High TEB             &                 \\ \hline
\multicolumn{1}{c|}{}                & Low tax        & \multicolumn{1}{l|}{}           & Low tax        & \multicolumn{1}{l|}{}          & Low tax             &                 \\
\multicolumn{1}{l|}{Smallest effect} & Low incentive  & \multicolumn{1}{l|}{+ 28.00\%} & Low incentive  & \multicolumn{1}{l|}{+ 58.69 \%} & Low incentive       & + 99.94 \%     \\
\multicolumn{1}{l|}{}                & Low TEB        & \multicolumn{1}{l|}{}           & Low TEB        & \multicolumn{1}{l|}{}          & Low TEB             &                 \\ \hline \hline
\multicolumn{7}{c}{High GEB}                                                                                                                                                      \\ \hline \hline
\multicolumn{1}{l|}{}                & \multicolumn{2}{c|}{VRES}                        & \multicolumn{2}{c|}{Welfare}                    & \multicolumn{2}{c}{Generation amount} \\ \hline
\multicolumn{1}{c|}{}                & Low tax        & \multicolumn{1}{l|}{}           & Low tax        & \multicolumn{1}{l|}{}          & Low tax             &                 \\
\multicolumn{1}{l|}{Largest effect}  & High incentive & \multicolumn{1}{l|}{+ 56.04 \%} & High incentive & \multicolumn{1}{l|}{+ 78.32 \%} & High incentive      & + 165.24 \%     \\
\multicolumn{1}{l|}{}                & High TEB       & \multicolumn{1}{l|}{}           & High TEB       & \multicolumn{1}{l|}{}          & High TEB            &                 \\ \hline
\multicolumn{1}{c|}{}                & Low tax        & \multicolumn{1}{l|}{}           & Low tax        & \multicolumn{1}{l|}{}          & Low tax             &                 \\
\multicolumn{1}{l|}{Smallest effect} & Low incentive  & \multicolumn{1}{l|}{+ 53.32 \%} & Low incentive  & \multicolumn{1}{l|}{+ 71.05 \%} & Low incentive       & +147.43 \%      \\
\multicolumn{1}{l|}{}                & Low TEB        & \multicolumn{1}{l|}{}           & Low TEB        & \multicolumn{1}{l|}{}          & Low TEB             &                 \\ \hline
\end{tabular}
\caption{Nordics case study: summary. Percentage values (\%) are calculated relative to the baseline case.}
\label{tab: nordic_summary}
\end{table}

It is important to bear in mind that the Nordic energy system has a low carbon footprint in 2018--2020. Nordic Energy Research \citep{oslo_nordics_nodate} indicates that the Nordics have reached about a 55.00\% renewable energy share in the generation mix in 2019. As the numerical results suggest, the optimal generation mix for combined Nordic and Baltic countries has about 35.84\% share of VRE when GEB and TEB are set to zero, and no carbon tax nor incentive are involved (baseline case). The results in Figure \ref{fig: nordics_tax_vres} demonstrate that using the ``worst'' composition of the input parameters values, i.e., ``low'' values for incentive, carbon tax and TEB, one can increase the share of VRE by about 28.0\% (that is, from 35.84\% in the baseline case to 63.84\%) compared to the case with no policies involved for the GenCos with ``low'' GEB and by about 53.3\% (i.e., 89.16 - 35.84 = 53.32\%) for the GenCos with ``high'' GEB. At the same time, Figures \ref{fig: nordics_incentive_welfare} and \ref{fig: nordics_incentive_generation} indicate that with the ``low'' values of the input parameters, the total generation and welfare values increase by roughly 99.9\% and 58.7\%, respectively, for the GenCos with a ``small'' GEB and by roughly 147.4\% and 71.1\%, accordingly for the GenCos with ``high'' GEB, compared to the baseline case. 

Figures \ref{fig: nordics_tax_vres}, \ref{fig: nordics_incentive_welfare} and \ref{fig: nordics_incentive_generation} also suggest that further increasing the carbon tax does not lead to any meaningful improvement regardless of the GEB value that GenCos possess. Potentially, this is because the VRE share in the generation mix has already exceeded 50\%. Hence, even if the expenses for the remaining share of the conventional generation rise, its impact on Gencos' profit is not crucial enough to motivate GenCos to change their investment decisions.

One can also consider the increase of the GEB as a renewables-driven policy. As Figures \ref{fig: nordics_tax_vres}, and \ref{fig: nordics_incentive_generation} suggest, by only increasing the GEB value from ``low'' to ``high'' while considering the rest of the input parameters (TEB, incentive and carbon tax) at a ``low'' values, one can increase (compared to the baseline case) the VRE share in the total generation mix by about 25.3\% (89.16 - 63.84 = 25.32\%) and the total generation amount by roughly 47.5\% (147.43 - 99.94 = 47.49\%) compared to the case of ``low'' GEB. However, the effect of such an approach on the total welfare is more modest as Figure \ref{fig: nordics_incentive_welfare} reports an increase of welfare (compared to the baseline) that is only about 12.4\% (71.05 - 58.69 = 12.36\%) higher than in the case with ``low'' GEB. The latter is due to the linear dependence between the demand and energy price, leading to the higher generation amount in case of higher GEB value and lower market price. 

% that a higher GEB value on its own leads to an increase in VRE share by
% A similar situation occurs with regard to the GEB. Considering the ``low'' value for GEB as in Table \ref{Tab: Nordic_case_input_parameters} does not lead to GEB-related constraints being active for any GenCo. Hence, considering a ``high'' value for GEB only leads to a slightly different GenCos investment decision portfolio. At the same time, the total investment expenses add up to the same value as in the case when the GEB value was ``low'' and hence, do not impact any of the output factors. As an example, Table \ref{tab: nordic case study investment portfolio} demonstrates the lack of differences in accumulated optimal investments of GenCos when considering ``high'' values for TEB, carbon tax and incentive but different values of GEB.

% \begin{table}[H!]
% \centering
% \begin{tabular}{l|cc}
%         & ``low'' GEB   & ``high'' GEB  \\  \hline \hline
% GenCo 1 & \euro2.9M & \euro3.1M \\
% GenCo 2 & \euro3.0M & \euro3.1M  \\
% GenCo 3 & \euro3.1M  & \euro3.1M  \\
% GenCo 4 & \euro2.8M & \euro2.6M \\
% GenCo 5 & \euro3.0M & \euro2.9M \\
% Total   & \euro14.8M & \euro14.8M  \\ \hline \hline
% \end{tabular}
% \caption{Optimal investments considering ``high'' values for TEB, carbon tax and incentive.}
% \label{tab: nordic case study investment portfolio}
% \end{table}

The two remaining input parameters that have an impact on the output factors under certain circumstances are the TEB and the VRE incentive. We first concentrate on the case when GEB is ``low''. Figures \ref{fig: nordics_tax_vres}, \ref{fig: nordics_incentive_welfare} and \ref{fig: nordics_incentive_generation} suggest that doubling the ``low'' TEB value does not impact any of the output factors if applied alone while the rest of the input parameters remaining at ``low'' values. Meanwhile, Figure \ref{fig: nordics_tax_vres} reports an increase (compared to the baseline) of VRE share in total generation mix by about 13.0\% (76.80 - 63.84 = 12.96\%) more than in the case with all ``low'' input parameters if one was to only increase the VRE incentive value from ``low'' to ``high''. Nevertheless, as Figures \ref{fig: nordics_incentive_welfare} and \ref{fig: nordics_incentive_generation} suggest, the impact of solely increasing the incentive value from ``low'' to ``high'' has a minor effect on the welfare and total generation values. The numerical experiments suggest that the increase (compared to the baseline) in total welfare and generation amount values is only higher by 7.9\% (66.58 - 58.69 = 7.89\%) and 5.0\% (104.96 - 99.94 = 5.02\%), respectively than in the case with all ``low'' input parameters values. Additionally, Figures \ref{fig: nordics_tax_vres}, \ref{fig: nordics_incentive_welfare} and \ref{fig: nordics_incentive_generation} suggest that considering the simultaneous increase of incentive and TEB values brings up an increase (compared to the baseline) of VRE share, total generation amount and total welfare that is higher by about 14.3\% (78.12 - 63.84 = 14.28\%), 5.5\% (105.41 - 99.94 = 5.47\%) and 8.1\% (66.75 - 58.69 = 8.06\%), respectively, than in the case when all the input parameters are at the ``low'' values.

In the case of the GenCos possessing ``high'' GEB, the inference regarding the effect of TEB and incentive on the output factor slightly changes. First of all, Figures \ref{fig: nordics_tax_vres}, \ref{fig: nordics_incentive_welfare} and \ref{fig: nordics_incentive_generation} suggest that in case of ``high'' GEB, solely increasing TEB value while keeping the incentive and the tax at the ``low'' values impacts all of the output factors. The increases (compared to the baseline) in VRE share, total generation amount and total welfare are higher by about 1.6\% (90.80 - 89.16 = 1.64\%), 9.9\% (157.31 - 147.43 = 9.88\%) and roughly by 1.0\% (72.00 - 71.05 = 0.95\%), respectively, than in the case with all (but the GEB) input parameters being at the ``low'' values. Increasing only the incentive value from ``low'' to ``high'' in the presence of the ``high'' GEB value and the rest of the input parameters being at the ``low'' values suggests a similar effect as in the case of solely increasing the ''TEB'' value. Figures \ref{fig: nordics_tax_vres}, \ref{fig: nordics_incentive_welfare} and \ref{fig: nordics_incentive_generation} suggest the increase (compared to the baseline) of the VRE share, total generation and total welfare is higher by about 1.3\% (90.44 - 89.16 = 1.28\%), 6.7\% (154.09 - 147.43 = 6.66\%) and 6.2\% (77.26 - 71.05 = 6.21\%), respectively, compared to the case when all the parameters (but the GEB) are ``low''. The simultaneous increase of TEB and incentive values considering ``high'' GEB value leads to the strongest increase in all the output parameters when compared to the baseline. Figures \ref{fig: nordics_tax_vres}, \ref{fig: nordics_incentive_welfare} and \ref{fig: nordics_incentive_generation} report the increase (compared to the baseline) of VRE share, welfare, and total generation is higher by about 2.7\% (91.88 - 89.16 = 2.72\%), 7.3\% (78.32 - 71.05 = 7.27\%) and 17.8\% (165.24 - 147.43 = 17.81\%), respectively than in the case with all (but the GEB) parameters being ``low''. 

\subsection{Summary of the results}

To conclude, the numerical results have shown that the ``low'' value of the carbon tax presented in Table \ref{Tab: Nordic_case_input_parameters} is efficient enough as the further increase of its value does not lead to an improvement in any of the output factors. Nevertheless, in case GenCos possess a ``low'' GEB, considering ``low'' values for TEB and an incentive leads to an increase by about 28.0\% in the VRE share in the optimal generation mix when compared with the baseline. At the same time, one can notice an increase by about 58.7\% in total welfare and an increase by roughly 99.94\% in total generation amount. If one is aimed at maximising VRE share, welfare and total generation values, one should consider simultaneously increasing TEB and incentive values to ``high''. Such an approach leads to an increase (compared to the baseline) of VRE share, welfare and total generation values, which is higher by about 14.3\%, 8.1\% and 5.5\%, respectively, than in the case if all the input parameters are ``low''. However, one should bear in mind that the numerical results suggest that ``high'' incentive value has a greater influence on VRE share and welfare increase if applied solely compared to only increasing TEB value.

In case GenCos possess a ``high'' GEB exploring ``low'' values of all the input parameters leads to an increase by about 53.3\%, 71.1\% and 147.4\% in VRE share, welfare and generation amount, respectively, when compared to the baseline results. Further increasing TEB and incentive to ``high'' values allows one to obtain the highest values for all the output factors. The results show an increase (compared to the baseline) of VRE share, welfare and generation amount that is higher by about 2.7\%, 7.3\% and 17.8\%, respectively, than in the case with all (but the GEB) input parameters being at the ``low'' value. Meanwhile, the results also suggest that if one was to solely increase just either TEB or incentive value, then the choice should be made in favour of TEB in case the aim is to increase the total generation amount or VRE share in the total generation mix. Alternatively, only increasing incentive value would have a greater impact on welfare than solely increasing TEB value when compared to the optimal welfare value for the case with all (but the GEB) parameters being at a ``low'' value.

\section{Conclusions}
\label{sec: conclusion}
In this paper, we study the impact of the TSO infrastructure expansion decisions in combination with carbon taxes and renewable-driven investment incentives on the optimal generation mix. To examine the impact of renewables-driven policies we propose a novel bi-level modelling assessment to plan optimal transmission infrastructure expansion. At the lower level, we consider a perfectly competitive energy market comprising GenCos who decide optimal generation levels and their own infrastructure expansion strategy. The upper level consists of a TSO who proactively anticipates the aforementioned decisions and decides the optimal transmission capacity expansion plan. To supplement the TSO decisions with other renewable-driven policies, we introduced carbon taxes and renewable capacity investment incentives in the model. Additionally, we accounted for variations in GenCos' and TSO's willingness to expand the infrastructure by introducing an upper limit on the generation (GEB) and transmission capacity expansion (TEB) costs. Therefore, as the input parameters for the proposed bi-level model, we considered different values of TEB, GEB, incentives and carbon tax. This paper examined the proposed modelling approach by applying it to a simple, three-node illustrative case study and a more realistic energy system representing Nordic and Baltic countries. The output factors explored in the analysis are the optimal total welfare, the share of VRE in the optimal generation mix and the total amount of energy generated.

In the illustrative case study, each of the parameters has been individually tested on a wide range of discrete values. The numerical experiments revealed the inefficiency of attempting to distinctly increase the budget allocated for transmission infrastructure expansion in case GenCos do not plan significant generation capacity expansions. Moreover, even when the generation expansion budget is large, the trade-off between the increase of VRE share and transmission infrastructure investment costs still poses the question of whether such a policy is efficient. Furthermore, the experiments identified difficulties in qualifying the impact of increasing the VRE capacity expansion incentives on each of the output factors due to the inconsistencies in the output data. However, the increase in carbon taxes has demonstrated a straightforward connection with the increase in VRE share once the marginal value of conventional energy becomes less profitable for GenCos. 

The numerical experiments for a more realistic case study depicting Nordic and Baltic countries' energy systems indicated the largest increase in all the output factors when compared to the case with no policies involved (i.e., when GEB, TEB, carbon tax and incentive are zeros) if one was to consider ``high'' values of TEB and incentive and ``low'' value for carbon tax regardless of the GEB value. However, the extent of increase in the values of all output factors followed by changing values of TEB and incentive from ``low'' to ``high'' varies depending on the GEB value. Therefore, the results indicate the prominent role of the TSO in the renewable transition of the energy market.

It is worth mentioning that to ensure the computational tractability of the proposed model, we allowed for a number of simplifications in the modelling approach. In particular, our model does not account for an intranodal spatial resolution. Therefore, the estimates derived from the proposed modelling methodology should be considered as an initial approximation, and further refinement may be necessary for practical implementation. In particular, the simplifications we considered are i) considering only Kirchhoff's voltage law, ii) lack of account for revenue streams for the TSO and iii) simplification of planning horizon and truncation of the scale of uncertainty representation. Nevertheless, the aforementioned simplifications do not diminish policy-related insights concluded from the numerical experiments. In addition, applying the proposed modelling assessment to the case study representing Nordic and Baltic countries' energy systems revealed the limitations associated with solving such instances using state-of-the-art solvers. Therefore, we discretised and enumerated possible upper-level decisions and afterwards, among all the possible alternatives, they determined the investment portfolio maximising the welfare. For the input parameters, we considered two types of values named ``low'' and ``high'' in accordance with open data regarding Nordic GenCos and TSO's infrastructure expansion investments plans and current EU member's carbon taxation systems.

Regarding further research, one could pinpoint a few possible directions. The first one is related to considering an imperfectly competitive market defining a Cornout oligopoly \citep{ruffin1971cournot} instead of perfect competition. The imperfectly competitive market structure may more closely represent the reality \citep{oikonomou2009white} as the oligopolies have access to information helping them in the decision-making process. However, this would require formulating a set of separate subproblems for each of the GenCos at the lower level which would further increase the computational burden. An alternative avenue for future research could involve the integration of storage batteries into the grid. This would potentially alleviate the intermittency problem regarding the availability of variable renewable energy sources and hence provide new insights, thereby offering novel perspectives on the outcomes of numerical experiments. However, one should bear in mind that this would also most probably further complicate the computational tractability of the problem. Another possible enhancement could stem from the development of an efficient solution method that allows one to consider a continuous range of investment decisions for the TSO at the upper level. Lastly, if the modelling simplifications caused by limitations of state-of-the-art solvers are overcome one could investigate how policy-related insights differ in case the model allows for a higher level of detail for energy system representation, investment projects, uncertainty formulation and the planning horizon. 

%TC:ignore
\section{Acknowledgements}
Nikita Belyak and Fabricio Oliveira were supported by the Research Council of Finland (decision number 348092). Steven A. Gabriel was supported by a grant from the Civil Infrastructure Systems program at the National Science Foundation (Award \#2113891).
Part of the research was developed in the Young Scientists Summer Program at the International Institute for Applied Systems Analysis, Laxenburg (Austria) with financial support from the Finnish National Member Organization. Additionally, part of the research was conducted during a research visit to the Chalmers University of Technology, Gothenburg (Sweden). This research also benefited from the presentations at the European Conference on Stochastic Optimization Computational Management Science 2022, Venice (Italy) and the EURO 2022 conference, Espoo (Finland). The authors also acknowledge the computational resources provided by the Aalto University (Finland) Science-IT project.

%TC:endignore

\bibliographystyle{elsarticle-num-names}
\bibliography{referencias.bib}

\begin{thebibliography}{58}
\expandafter\ifx\csname natexlab\endcsname\relax\def\natexlab#1{#1}\fi
\providecommand{\url}[1]{\texttt{#1}}
\providecommand{\href}[2]{#2}
\providecommand{\path}[1]{#1}
\providecommand{\DOIprefix}{doi:}
\providecommand{\ArXivprefix}{arXiv:}
\providecommand{\URLprefix}{URL: }
\providecommand{\Pubmedprefix}{pmid:}
\providecommand{\doi}[1]{\href{http://dx.doi.org/#1}{\path{#1}}}
\providecommand{\Pubmed}[1]{\href{pmid:#1}{\path{#1}}}
\providecommand{\bibinfo}[2]{#2}
\ifx\xfnm\relax \def\xfnm[#1]{\unskip,\space#1}\fi
%Type = Article
\bibitem[{Taghvaee et~al.(2023)Taghvaee, Nodehi, Arani, Jafari, and Shirazi}]{taghvaee2023sustainability}
\bibinfo{author}{V.~M. Taghvaee}, \bibinfo{author}{M.~Nodehi}, \bibinfo{author}{A.~A. Arani}, \bibinfo{author}{Y.~Jafari}, \bibinfo{author}{J.~K. Shirazi},
\newblock \bibinfo{title}{Sustainability spillover effects of social, environment and economy: mapping global sustainable development in a systematic analysis},
\newblock \bibinfo{journal}{Asia-Pacific Journal of Regional Science} \bibinfo{volume}{7} (\bibinfo{year}{2023}) \bibinfo{pages}{329--353}.
%Type = Article
\bibitem[{Nasrollahi et~al.(2020)Nasrollahi, Hashemi, Bameri, and Mohamad~Taghvaee}]{nasrollahi2020environmental}
\bibinfo{author}{Z.~Nasrollahi}, \bibinfo{author}{M.-s. Hashemi}, \bibinfo{author}{S.~Bameri}, \bibinfo{author}{V.~Mohamad~Taghvaee},
\newblock \bibinfo{title}{Environmental pollution, economic growth, population, industrialization, and technology in weak and strong sustainability: using stirpat model},
\newblock \bibinfo{journal}{Environment, Development and Sustainability} \bibinfo{volume}{22} (\bibinfo{year}{2020}) \bibinfo{pages}{1105--1122}.
%Type = Article
\bibitem[{Taghvaee et~al.(2022)Taghvaee, Arani, Soretz, and Agheli}]{taghvaee2022comparing}
\bibinfo{author}{V.~M. Taghvaee}, \bibinfo{author}{A.~A. Arani}, \bibinfo{author}{S.~Soretz}, \bibinfo{author}{L.~Agheli},
\newblock \bibinfo{title}{Comparing energy efficiency and price policy from a sustainable development perspective: Using fossil fuel demand elasticities in iran},
\newblock \bibinfo{journal}{MRS Energy \& Sustainability} \bibinfo{volume}{9} (\bibinfo{year}{2022}) \bibinfo{pages}{480--493}.
%Type = Article
\bibitem[{Wolf et~al.(2021)Wolf, Teitge, Mielke, Sch{\"u}tze, and Jaeger}]{wolf_european_2021}
\bibinfo{author}{S.~Wolf}, \bibinfo{author}{J.~Teitge}, \bibinfo{author}{J.~Mielke}, \bibinfo{author}{F.~Sch{\"u}tze}, \bibinfo{author}{C.~Jaeger},
\newblock \bibinfo{title}{The {European} {Green} {Deal} --- {More} {Than} {Climate} {Neutrality}},
\newblock \bibinfo{journal}{Intereconomics} \bibinfo{volume}{56} (\bibinfo{year}{2021}) \bibinfo{pages}{99--107}. \URLprefix \url{https://doi.org/10.1007/s10272-021-0963-z}. \DOIprefix\doi{10.1007/s10272-021-0963-z}.
%Type = Misc
\bibitem[{{European Commission}(2021)}]{noauthor_climate_nodate}
\bibinfo{author}{{European Commission}}, \bibinfo{title}{Climate strategies \& targets}, \bibinfo{year}{2021}. \URLprefix \url{https://ec.europa.eu/clima/eu-action/climate-strategies-targets_en}.
%Type = Article
\bibitem[{Steffen(2018)}]{STEFFEN2018280}
\bibinfo{author}{B.~Steffen},
\newblock \bibinfo{title}{The importance of project finance for renewable energy projects},
\newblock \bibinfo{journal}{Energy Economics} \bibinfo{volume}{69} (\bibinfo{year}{2018}) \bibinfo{pages}{280--294}. \URLprefix \url{https://www.sciencedirect.com/science/article/pii/S0140988317303870}. \DOIprefix\doi{https://doi.org/10.1016/j.eneco.2017.11.006}.
%Type = Misc
\bibitem[{{Intergovernmental Panel on Climate Change}(2014)}]{Cambridge_climate}
\bibinfo{author}{{Intergovernmental Panel on Climate Change}}, \bibinfo{title}{Climate {Change} 2013 -- {The} {Physical} {Science} {Basis}}, \bibinfo{year}{2014}.
%Type = Misc
\bibitem[{Redl et~al.(2021)Redl, Hein, Buck, Graichen, and (Ember)}]{agora_european_2021}
\bibinfo{author}{C.~Redl}, \bibinfo{author}{F.~Hein}, \bibinfo{author}{M.~Buck}, \bibinfo{author}{D.~P. Graichen}, \bibinfo{author}{D.~J. (Ember)}, \bibinfo{title}{The european power sector in 2020: Up-to-date analysis on the electricity transition}, \bibinfo{year}{2021}.
%Type = Article
\bibitem[{Haar(2020)}]{HAAR2020111483}
\bibinfo{author}{L.~Haar},
\newblock \bibinfo{title}{An empirical analysis of the fiscal incidence of renewable energy support in the european union},
\newblock \bibinfo{journal}{Energy Policy} \bibinfo{volume}{143} (\bibinfo{year}{2020}) \bibinfo{pages}{111483}. \URLprefix \url{https://www.sciencedirect.com/science/article/pii/S0301421520302317}. \DOIprefix\doi{https://doi.org/10.1016/j.enpol.2020.111483}.
%Type = Article
\bibitem[{Sinsel et~al.(2020)Sinsel, Riemke, and Hoffmann}]{SINSEL20202271}
\bibinfo{author}{S.~R. Sinsel}, \bibinfo{author}{R.~L. Riemke}, \bibinfo{author}{V.~H. Hoffmann},
\newblock \bibinfo{title}{Challenges and solution technologies for the integration of variable renewable energy sources---a review},
\newblock \bibinfo{journal}{Renewable Energy} \bibinfo{volume}{145} (\bibinfo{year}{2020}) \bibinfo{pages}{2271--2285}. \URLprefix \url{https://www.sciencedirect.com/science/article/pii/S0960148119309875}. \DOIprefix\doi{https://doi.org/10.1016/j.renene.2019.06.147}.
%Type = Article
\bibitem[{Moreira et~al.(2017)Moreira, Pozo, Street, and Sauma}]{7752917}
\bibinfo{author}{A.~Moreira}, \bibinfo{author}{D.~Pozo}, \bibinfo{author}{A.~Street}, \bibinfo{author}{E.~Sauma},
\newblock \bibinfo{title}{Reliable renewable generation and transmission expansion planning: Co-optimizing system's resources for meeting renewable targets},
\newblock \bibinfo{journal}{IEEE Transactions on Power Systems} \bibinfo{volume}{32} (\bibinfo{year}{2017}) \bibinfo{pages}{3246--3257}. \DOIprefix\doi{10.1109/TPWRS.2016.2631450}.
%Type = Article
\bibitem[{Hemmati et~al.(2013)Hemmati, Hooshmand, and Khodabakhshian}]{HEMMATI2013312}
\bibinfo{author}{R.~Hemmati}, \bibinfo{author}{R.-A. Hooshmand}, \bibinfo{author}{A.~Khodabakhshian},
\newblock \bibinfo{title}{State-of-the-art of transmission expansion planning: Comprehensive review},
\newblock \bibinfo{journal}{Renewable and Sustainable Energy Reviews} \bibinfo{volume}{23} (\bibinfo{year}{2013}) \bibinfo{pages}{312--319}. \URLprefix \url{https://www.sciencedirect.com/science/article/pii/S1364032113001743}. \DOIprefix\doi{https://doi.org/10.1016/j.rser.2013.03.015}.
%Type = Article
\bibitem[{Lumbreras and Ramos(2016)}]{LUMBRERAS201619}
\bibinfo{author}{S.~Lumbreras}, \bibinfo{author}{A.~Ramos},
\newblock \bibinfo{title}{The new challenges to transmission expansion planning. survey of recent practice and literature review},
\newblock \bibinfo{journal}{Electric Power Systems Research} \bibinfo{volume}{134} (\bibinfo{year}{2016}) \bibinfo{pages}{19--29}. \URLprefix \url{https://www.sciencedirect.com/science/article/pii/S0378779615003090}. \DOIprefix\doi{https://doi.org/10.1016/j.epsr.2015.10.013}.
%Type = Inproceedings
\bibitem[{Niharika et~al.(2016)Niharika, Verma, and Mukherjee}]{7583779}
\bibinfo{author}{Niharika}, \bibinfo{author}{S.~Verma}, \bibinfo{author}{V.~Mukherjee},
\newblock \bibinfo{title}{Transmission expansion planning: A review},
\newblock in: \bibinfo{booktitle}{2016 International Conference on Energy Efficient Technologies for Sustainability (ICEETS)}, \bibinfo{year}{2016}, pp. \bibinfo{pages}{350--355}. \DOIprefix\doi{10.1109/ICEETS.2016.7583779}.
%Type = Article
\bibitem[{Sun et~al.(2018)Sun, Cremer, and Strbac}]{SUN2018546}
\bibinfo{author}{M.~Sun}, \bibinfo{author}{J.~Cremer}, \bibinfo{author}{G.~Strbac},
\newblock \bibinfo{title}{A novel data-driven scenario generation framework for transmission expansion planning with high renewable energy penetration},
\newblock \bibinfo{journal}{Applied Energy} \bibinfo{volume}{228} (\bibinfo{year}{2018}) \bibinfo{pages}{546--555}. \URLprefix \url{https://www.sciencedirect.com/science/article/pii/S0306261918309656}. \DOIprefix\doi{https://doi.org/10.1016/j.apenergy.2018.06.095}.
%Type = Article
\bibitem[{Mortaz and Valenzuela(2019)}]{MORTAZ201935}
\bibinfo{author}{E.~Mortaz}, \bibinfo{author}{J.~Valenzuela},
\newblock \bibinfo{title}{Evaluating the impact of renewable generation on transmission expansion planning},
\newblock \bibinfo{journal}{Electric Power Systems Research} \bibinfo{volume}{169} (\bibinfo{year}{2019}) \bibinfo{pages}{35--44}. \URLprefix \url{https://www.sciencedirect.com/science/article/pii/S0378779618304048}. \DOIprefix\doi{https://doi.org/10.1016/j.epsr.2018.12.007}.
%Type = Article
\bibitem[{Tian et~al.(2020)Tian, Sun, Han, and Yang}]{8945251}
\bibinfo{author}{K.~Tian}, \bibinfo{author}{W.~Sun}, \bibinfo{author}{D.~Han}, \bibinfo{author}{C.~Yang},
\newblock \bibinfo{title}{Coordinated planning with predetermined renewable energy generation targets using extended two-stage robust optimization},
\newblock \bibinfo{journal}{IEEE Access} \bibinfo{volume}{8} (\bibinfo{year}{2020}) \bibinfo{pages}{2395--2407}. \DOIprefix\doi{10.1109/ACCESS.2019.2962841}.
%Type = Article
\bibitem[{Zhang et~al.(2020)Zhang, Cheng, Liu, Zhang, Zhang, and Li}]{ZHANG2020105944}
\bibinfo{author}{C.~Zhang}, \bibinfo{author}{H.~Cheng}, \bibinfo{author}{L.~Liu}, \bibinfo{author}{H.~Zhang}, \bibinfo{author}{X.~Zhang}, \bibinfo{author}{G.~Li},
\newblock \bibinfo{title}{Coordination planning of wind farm, energy storage and transmission network with high-penetration renewable energy},
\newblock \bibinfo{journal}{International Journal of Electrical Power \& Energy Systems} \bibinfo{volume}{120} (\bibinfo{year}{2020}) \bibinfo{pages}{105944}. \URLprefix \url{https://www.sciencedirect.com/science/article/pii/S0142061519324445}. \DOIprefix\doi{https://doi.org/10.1016/j.ijepes.2020.105944}.
%Type = Article
\bibitem[{Tatar et~al.(2024)Tatar, Harati, Farokhi, Taghvaee, and Wilson}]{tatar2024good}
\bibinfo{author}{M.~Tatar}, \bibinfo{author}{J.~Harati}, \bibinfo{author}{S.~Farokhi}, \bibinfo{author}{V.~Taghvaee}, \bibinfo{author}{F.~A. Wilson},
\newblock \bibinfo{title}{Good governance and natural resource management in oil and gas resource-rich countries: A machine learning approach},
\newblock \bibinfo{journal}{Resources Policy} \bibinfo{volume}{89} (\bibinfo{year}{2024}) \bibinfo{pages}{104583}.
%Type = Article
\bibitem[{Zhang et~al.(2016)Zhang, Hu, Springer, Li, and Shen}]{ZHANG201684}
\bibinfo{author}{N.~Zhang}, \bibinfo{author}{Z.~Hu}, \bibinfo{author}{C.~Springer}, \bibinfo{author}{Y.~Li}, \bibinfo{author}{B.~Shen},
\newblock \bibinfo{title}{A bi-level integrated generation-transmission planning model incorporating the impacts of demand response by operation simulation},
\newblock \bibinfo{journal}{Energy Conversion and Management} \bibinfo{volume}{123} (\bibinfo{year}{2016}) \bibinfo{pages}{84--94}. \URLprefix \url{https://www.sciencedirect.com/science/article/pii/S019689041630499X}. \DOIprefix\doi{https://doi.org/10.1016/j.enconman.2016.06.020}.
%Type = Article
\bibitem[{Virasjoki et~al.(2020)Virasjoki, Siddiqui, Oliveira, and Salo}]{VIRASJOKI2020104716}
\bibinfo{author}{V.~Virasjoki}, \bibinfo{author}{A.~S. Siddiqui}, \bibinfo{author}{F.~Oliveira}, \bibinfo{author}{A.~Salo},
\newblock \bibinfo{title}{Utility-scale energy storage in an imperfectly competitive power sector},
\newblock \bibinfo{journal}{Energy Economics} \bibinfo{volume}{88} (\bibinfo{year}{2020}) \bibinfo{pages}{104716}.
%Type = Article
\bibitem[{Siddiqui et~al.(2019)Siddiqui, Tanaka, and Chen}]{SIDDIQUI2019208}
\bibinfo{author}{A.~S. Siddiqui}, \bibinfo{author}{M.~Tanaka}, \bibinfo{author}{Y.~Chen},
\newblock \bibinfo{title}{Sustainable transmission planning in imperfectly competitive electricity industries: Balancing economic and environmental outcomes},
\newblock \bibinfo{journal}{European Journal of Operational Research} \bibinfo{volume}{275} (\bibinfo{year}{2019}) \bibinfo{pages}{208--223}. \URLprefix \url{https://www.sciencedirect.com/science/article/pii/S0377221718309573}. \DOIprefix\doi{https://doi.org/10.1016/j.ejor.2018.11.032}.
%Type = Article
\bibitem[{Quintela et~al.(2009)Quintela, Redondo, Melchor, and Redondo}]{4909474}
\bibinfo{author}{F.~R. Quintela}, \bibinfo{author}{R.~C. Redondo}, \bibinfo{author}{N.~R. Melchor}, \bibinfo{author}{M.~Redondo},
\newblock \bibinfo{title}{A general approach to kirchhoff's laws},
\newblock \bibinfo{journal}{IEEE Transactions on Education} \bibinfo{volume}{52} (\bibinfo{year}{2009}) \bibinfo{pages}{273--278}. \DOIprefix\doi{10.1109/TE.2008.928189}.
%Type = Article
\bibitem[{Padiyar and Shanbhag(1988)}]{PADIYAR198817}
\bibinfo{author}{K.~Padiyar}, \bibinfo{author}{R.~Shanbhag},
\newblock \bibinfo{title}{Comparison of methods for transmission system expansion using network flow and dc load flow models},
\newblock \bibinfo{journal}{International Journal of Electrical Power \& Energy Systems} \bibinfo{volume}{10} (\bibinfo{year}{1988}) \bibinfo{pages}{17--24}. \URLprefix \url{https://www.sciencedirect.com/science/article/pii/014206158890004X}. \DOIprefix\doi{https://doi.org/10.1016/0142-0615(88)90004-X}.
%Type = Book
\bibitem[{Gabriel et~al.(2013)Gabriel, Conejo, Fuller, Hobbs, and Ruiz}]{gabriel_complementarity_2013}
\bibinfo{author}{S.~A. Gabriel}, \bibinfo{author}{A.~J. Conejo}, \bibinfo{author}{J.~D. Fuller}, \bibinfo{author}{B.~F. Hobbs}, \bibinfo{author}{C.~Ruiz}, \bibinfo{title}{Complementarity {Modeling} in {Energy} {Markets}}, International {Series} in {Operations} {Research} \& {Management} {Science}, \bibinfo{publisher}{Springer-Verlag}, \bibinfo{address}{New York}, \bibinfo{year}{2013}. \URLprefix \url{https://www.springer.com/gp/book/9781441961228}. \DOIprefix\doi{10.1007/978-1-4419-6123-5}.
%Type = Article
\bibitem[{Bezanson et~al.(2017)Bezanson, Edelman, Karpinski, and Shah}]{Bazanson2017Julia}
\bibinfo{author}{J.~Bezanson}, \bibinfo{author}{A.~Edelman}, \bibinfo{author}{S.~Karpinski}, \bibinfo{author}{V.~B. Shah},
\newblock \bibinfo{title}{Julia: {A} {Fresh} {Approach} to {Numerical} {Computing}},
\newblock \bibinfo{journal}{SIAM Review} \bibinfo{volume}{59} (\bibinfo{year}{2017}) \bibinfo{pages}{65--98}. \URLprefix \url{https://epubs.siam.org/doi/10.1137/141000671}. \DOIprefix\doi{10.1137/141000671}.
%Type = Misc
\bibitem[{Garcia et~al.(2022)Garcia, Bodin, and Street}]{Joaquim_bilevel_2022}
\bibinfo{author}{J.~D. Garcia}, \bibinfo{author}{G.~Bodin}, \bibinfo{author}{A.~Street}, \bibinfo{title}{Bileveljump.jl: Modeling and solving bilevel optimization in julia}, \bibinfo{year}{2022}. \URLprefix \url{https://arxiv.org/abs/2205.02307}. \DOIprefix\doi{10.48550/ARXIV.2205.02307}.
%Type = Misc
\bibitem[{Gurobi~Optimization(2020)}]{gurobi}
\bibinfo{author}{L.~Gurobi~Optimization}, \bibinfo{title}{Gurobi optimizer reference manual}, \bibinfo{year}{2020}. \URLprefix \url{http://www.gurobi.com}.
%Type = Misc
\bibitem[{Belyak(2022{\natexlab{a}})}]{belyak_ill_2022}
\bibinfo{author}{N.~Belyak}, \bibinfo{title}{Transmission system expansion planning: illustrative example}, \bibinfo{howpublished}{\url{https://github.com/gamma-opt/TSEP_illustrative_example}}, \bibinfo{year}{2022}{\natexlab{a}}.
%Type = Misc
\bibitem[{Belyak(2022{\natexlab{b}})}]{belyak_nor_2022}
\bibinfo{author}{N.~Belyak}, \bibinfo{title}{Transmission system expansion planning: Nordic case study}, \bibinfo{howpublished}{\url{https://github.com/gamma-opt/TSEP_Nordic_case_study}}, \bibinfo{year}{2022}{\natexlab{b}}.
%Type = Article
\bibitem[{Mir et~al.(2020)Mir, Alghassab, Ullah, Khan, Lu, and Imran}]{mir2020review}
\bibinfo{author}{A.~A. Mir}, \bibinfo{author}{M.~Alghassab}, \bibinfo{author}{K.~Ullah}, \bibinfo{author}{Z.~A. Khan}, \bibinfo{author}{Y.~Lu}, \bibinfo{author}{M.~Imran},
\newblock \bibinfo{title}{A review of electricity demand forecasting in low and middle income countries: The demand determinants and horizons},
\newblock \bibinfo{journal}{Sustainability} \bibinfo{volume}{12} (\bibinfo{year}{2020}) \bibinfo{pages}{5931}.
%Type = Incollection
\bibitem[{Timilsina and Shah(2022)}]{timilsina2022economics}
\bibinfo{author}{G.~R. Timilsina}, \bibinfo{author}{K.~U. Shah},
\newblock \bibinfo{title}{Economics of renewable energy: A comparison of electricity production costs across technologies},
\newblock in: \bibinfo{booktitle}{Oxford Research Encyclopedia of Environmental Science}, \bibinfo{year}{2022}.
%Type = Misc
\bibitem[{{d-maps}(2022)}]{noauthor_scandinavia_nodate}
\bibinfo{author}{{d-maps}}, \bibinfo{title}{Scandinavia free map}, \bibinfo{year}{2022}. \URLprefix \url{https://d-maps.com/carte.php?num_car=5972&lang=en}.
%Type = Misc
\bibitem[{Fortum(2022)}]{fortum_country_nodate}
\bibinfo{author}{Fortum}, \bibinfo{title}{Fortum worldwide report: countries fact sheets}, \bibinfo{year}{2022}. \URLprefix \url{https://www.fortum.com/about-us/our-company/fortum-worldwide/country-fact-sheets}.
%Type = Misc
\bibitem[{IBM(2021)}]{noauthor_ibm_2021}
\bibinfo{author}{IBM}, \bibinfo{title}{{IBM} {Documentation}}, \bibinfo{year}{2021}. \URLprefix \url{https://www.ibm.com/docs/en/icos/12.10.0?topic=cplex-overview-apis}.
%Type = Misc
\bibitem[{GAMS(2023)}]{noauthor_gams_nodate}
\bibinfo{author}{GAMS}, \bibinfo{title}{{GAMS} {Documentation} {Center}}, \bibinfo{year}{2023}. \URLprefix \url{https://www.gams.com/latest/docs/}.
%Type = Article
\bibitem[{Li and Pye(2018)}]{LI2018965}
\bibinfo{author}{P.-H. Li}, \bibinfo{author}{S.~Pye},
\newblock \bibinfo{title}{Assessing the benefits of demand-side flexibility in residential and transport sectors from an integrated energy systems perspective},
\newblock \bibinfo{journal}{Applied Energy} \bibinfo{volume}{228} (\bibinfo{year}{2018}) \bibinfo{pages}{965--979}. \URLprefix \url{https://www.sciencedirect.com/science/article/pii/S0306261918310237}. \DOIprefix\doi{https://doi.org/10.1016/j.apenergy.2018.06.153}.
%Type = Article
\bibitem[{Lind and Espegren(2017)}]{LIND201744}
\bibinfo{author}{A.~Lind}, \bibinfo{author}{K.~Espegren},
\newblock \bibinfo{title}{The use of energy system models for analysing the transition to low-carbon cities – the case of oslo},
\newblock \bibinfo{journal}{Energy Strategy Reviews} \bibinfo{volume}{15} (\bibinfo{year}{2017}) \bibinfo{pages}{44--56}. \URLprefix \url{https://www.sciencedirect.com/science/article/pii/S2211467X17300019}. \DOIprefix\doi{https://doi.org/10.1016/j.esr.2017.01.001}.
%Type = Article
\bibitem[{Fodstad et~al.(2022)Fodstad, {Crespo del Granado}, Hellemo, Knudsen, Pisciella, Silvast, Bordin, Schmidt, and Straus}]{FODSTAD2022112246}
\bibinfo{author}{M.~Fodstad}, \bibinfo{author}{P.~{Crespo del Granado}}, \bibinfo{author}{L.~Hellemo}, \bibinfo{author}{B.~R. Knudsen}, \bibinfo{author}{P.~Pisciella}, \bibinfo{author}{A.~Silvast}, \bibinfo{author}{C.~Bordin}, \bibinfo{author}{S.~Schmidt}, \bibinfo{author}{J.~Straus},
\newblock \bibinfo{title}{Next frontiers in energy system modelling: A review on challenges and the state of the art},
\newblock \bibinfo{journal}{Renewable and Sustainable Energy Reviews} \bibinfo{volume}{160} (\bibinfo{year}{2022}) \bibinfo{pages}{112246}. \URLprefix \url{https://www.sciencedirect.com/science/article/pii/S136403212200168X}. \DOIprefix\doi{https://doi.org/10.1016/j.rser.2022.112246}.
%Type = Article
\bibitem[{Teichgraeber et~al.(2019)Teichgraeber, Kuepper, and Brandt}]{Teichgraeber2019joss}
\bibinfo{author}{H.~Teichgraeber}, \bibinfo{author}{L.~E. Kuepper}, \bibinfo{author}{A.~R. Brandt},
\newblock \bibinfo{title}{Timeseriesclustering : An extensible framework in julia},
\newblock \bibinfo{journal}{Journal of Open Source Software} \bibinfo{volume}{4} (\bibinfo{year}{2019}) \bibinfo{pages}{1573}. \DOIprefix\doi{https://doi.org/10.21105/joss.01573}.
%Type = Misc
\bibitem[{{Nordic Energy Research}(2022)}]{oslo_nordics_nodate}
\bibinfo{author}{{Nordic Energy Research}}, \bibinfo{title}{Nordics lead {Europe} in renewables – {Nordic} {Energy} {Research}}, \bibinfo{year}{2022}. \URLprefix \url{https://www.nordicenergy.org/article/nordics-lead-europe-in-renewables/}.
%Type = Article
\bibitem[{Ruffin(1971)}]{ruffin1971cournot}
\bibinfo{author}{R.~J. Ruffin},
\newblock \bibinfo{title}{Cournot oligopoly and competitive behaviour},
\newblock \bibinfo{journal}{The Review of Economic Studies} \bibinfo{volume}{38} (\bibinfo{year}{1971}) \bibinfo{pages}{493--502}.
%Type = Inproceedings
\bibitem[{Oikonomou et~al.(2009)Oikonomou, Di~Giacomo, Russolillo, Becchis et~al.}]{oikonomou2009white}
\bibinfo{author}{V.~Oikonomou}, \bibinfo{author}{M.~Di~Giacomo}, \bibinfo{author}{D.~Russolillo}, \bibinfo{author}{F.~Becchis}, et~al.,
\newblock \bibinfo{title}{White certificates in an oligopoly market: closer to reality?},
\newblock in: \bibinfo{booktitle}{Proceedings of the ECEEE summer study}, \bibinfo{year}{2009}, pp. \bibinfo{pages}{1071--1080}.
%Type = Article
\bibitem[{Hirth et~al.(2018)Hirth, M{\"u}hlenpfordt, and Bulkeley}]{hirth_entso-e_2018}
\bibinfo{author}{L.~Hirth}, \bibinfo{author}{J.~M{\"u}hlenpfordt}, \bibinfo{author}{M.~Bulkeley},
\newblock \bibinfo{title}{The {ENTSO}-{E} {Transparency} {Platform} -- {A} review of {Europe}'s most ambitious electricity data platform},
\newblock \bibinfo{journal}{Applied Energy} \bibinfo{volume}{225} (\bibinfo{year}{2018}) \bibinfo{pages}{1054--1067}. \URLprefix \url{https://www.sciencedirect.com/science/article/pii/S0306261918306068}. \DOIprefix\doi{10.1016/j.apenergy.2018.04.048}.
%Type = Article
\bibitem[{John(1999)}]{john1999technical}
\bibinfo{author}{M.~John},
\newblock \bibinfo{title}{Technical analysis of the financial markets},
\newblock \bibinfo{journal}{New York Institute of Finance} \bibinfo{volume}{24} (\bibinfo{year}{1999}) \bibinfo{pages}{1--5}.
%Type = Article
\bibitem[{Kan et~al.(2020)Kan, Hedenus, and Reichenberg}]{KAN2020117015}
\bibinfo{author}{X.~Kan}, \bibinfo{author}{F.~Hedenus}, \bibinfo{author}{L.~Reichenberg},
\newblock \bibinfo{title}{The cost of a future low-carbon electricity system without nuclear power – the case of sweden},
\newblock \bibinfo{journal}{Energy} \bibinfo{volume}{195} (\bibinfo{year}{2020}) \bibinfo{pages}{117015}. \URLprefix \url{https://www.sciencedirect.com/science/article/pii/S0360544220301225}. \DOIprefix\doi{https://doi.org/10.1016/j.energy.2020.117015}.
%Type = Article
\bibitem[{Chatzimouratidis and Pilavachi(2009)}]{CHATZIMOURATIDIS2009778}
\bibinfo{author}{A.~I. Chatzimouratidis}, \bibinfo{author}{P.~A. Pilavachi},
\newblock \bibinfo{title}{Technological, economic and sustainability evaluation of power plants using the analytic hierarchy process},
\newblock \bibinfo{journal}{Energy Policy} \bibinfo{volume}{37} (\bibinfo{year}{2009}) \bibinfo{pages}{778--787}. \URLprefix \url{https://www.sciencedirect.com/science/article/pii/S0301421508005880}. \DOIprefix\doi{https://doi.org/10.1016/j.enpol.2008.10.009}.
%Type = Article
\bibitem[{Cui et~al.(2019)Cui, Hultman, Edwards, He, Sen, Surana, McJeon, Iyer, Patel, Yu, Nace, and Shearer}]{Cui_2019}
\bibinfo{author}{R.~Y. Cui}, \bibinfo{author}{N.~Hultman}, \bibinfo{author}{M.~R. Edwards}, \bibinfo{author}{L.~He}, \bibinfo{author}{A.~Sen}, \bibinfo{author}{K.~Surana}, \bibinfo{author}{H.~McJeon}, \bibinfo{author}{G.~Iyer}, \bibinfo{author}{P.~Patel}, \bibinfo{author}{S.~Yu}, \bibinfo{author}{T.~Nace}, \bibinfo{author}{C.~Shearer},
\newblock \bibinfo{title}{Quantifying operational lifetimes for coal power plants under the paris goals},
\newblock \bibinfo{journal}{Nature Communications} \bibinfo{volume}{10} (\bibinfo{year}{2019}) \bibinfo{pages}{4759}. \URLprefix \url{https://doi.org/10.1038/s41467-019-12618-3}. \DOIprefix\doi{10.1038/s41467-019-12618-3}.
%Type = Misc
\bibitem[{Rintamäki(2022)}]{Rintamaki_2022}
\bibinfo{author}{T.~Rintamäki}, \bibinfo{title}{Robust optimization for power markets}, \bibinfo{howpublished}{\url{https://github.com/tuomasr/robust-dev}}, \bibinfo{year}{2022}.
%Type = Misc
\bibitem[{Condeixa et~al.(2022)Condeixa, Oliveira, and Tollander~de Balsch}]{Condeixa_2022}
\bibinfo{author}{L.~Condeixa}, \bibinfo{author}{F.~Oliveira}, \bibinfo{author}{J.~Tollander~de Balsch}, \bibinfo{title}{Transmission capacity expansion problem modelling}, \bibinfo{howpublished}{\url{https://github.com/gamma-opt/EnergySystemModeling.jl}}, \bibinfo{year}{2022}.
%Type = Article
\bibitem[{Nycander et~al.(2020)Nycander, Söder, Olauson, and Eriksson}]{NYCANDER2020942}
\bibinfo{author}{E.~Nycander}, \bibinfo{author}{L.~Söder}, \bibinfo{author}{J.~Olauson}, \bibinfo{author}{R.~Eriksson},
\newblock \bibinfo{title}{Curtailment analysis for the nordic power system considering transmission capacity, inertia limits and generation flexibility},
\newblock \bibinfo{journal}{Renewable Energy} \bibinfo{volume}{152} (\bibinfo{year}{2020}) \bibinfo{pages}{942--960}. \URLprefix \url{https://www.sciencedirect.com/science/article/pii/S0960148120300641}. \DOIprefix\doi{https://doi.org/10.1016/j.renene.2020.01.059}.
%Type = Article
\bibitem[{Schlachtberger et~al.(2017)Schlachtberger, Brown, Schramm, and Greiner}]{SCHLACHTBERGER2017469}
\bibinfo{author}{D.~Schlachtberger}, \bibinfo{author}{T.~Brown}, \bibinfo{author}{S.~Schramm}, \bibinfo{author}{M.~Greiner},
\newblock \bibinfo{title}{The benefits of cooperation in a highly renewable european electricity network},
\newblock \bibinfo{journal}{Energy} \bibinfo{volume}{134} (\bibinfo{year}{2017}) \bibinfo{pages}{469--481}. \URLprefix \url{https://www.sciencedirect.com/science/article/pii/S0360544217309969}. \DOIprefix\doi{https://doi.org/10.1016/j.energy.2017.06.004}.
%Type = Misc
\bibitem[{{Nordic TSOs}(2022)}]{nordics_tsos}
\bibinfo{author}{{Nordic TSOs}}, \bibinfo{title}{2022 solutions for a green nordic energy system}, \bibinfo{year}{2022}. \URLprefix \url{chrome-extension://efaidnbmnnnibpcajpcglclefindmkaj/https://www.svk.se/siteassets/om-oss/rapporter/2022/solutions-report-2022.pdf}.
%Type = Misc
\bibitem[{Vattenfall(2022)}]{vattenfall_2022}
\bibinfo{author}{Vattenfall}, \bibinfo{title}{Our investment programme}, \bibinfo{year}{2022}. \URLprefix \url{https://group.vattenfall.com/investors/understanding-vattenfall/investment-plan}.
%Type = Misc
\bibitem[{{European court of auditors}(2022)}]{EU_review_taxation_2022}
\bibinfo{author}{{European court of auditors}}, \bibinfo{title}{Energy taxation, carbon pricing and energy subsidies 2022}, \bibinfo{year}{2022}. \URLprefix \url{https://www.eca.europa.eu/Lists/ECADocuments/RW22_01/RW_Energy_taxation_EN.pdf}.
%Type = Misc
\bibitem[{{Energy for Growth Hub}(2022)}]{noauthor_GT}
\bibinfo{author}{{Energy for Growth Hub}}, \bibinfo{title}{Should lower-income countries build open cycle or combined cycle gas turbines?}, \bibinfo{year}{2022}. \URLprefix \url{https://energyforgrowth.org/article/should-lower-income-countries-build-open-cycle-or-combined-cycle-gas-turbines/}.
%Type = Misc
\bibitem[{{U.S. Environmental Protection Agency}(2019)}]{noauthor_coal}
\bibinfo{author}{{U.S. Environmental Protection Agency}}, \bibinfo{title}{Greenhouse gas reporting program industrial profile: Power plants sector}, \bibinfo{year}{2019}. \URLprefix \url{https://www.epa.gov/sites/default/files/2020-12/documents/power_plants_2017_industrial_profile_updated_2020.pdf}.
%Type = Misc
\bibitem[{{Partnership for policy integrity}(2011)}]{noauthor_biomass}
\bibinfo{author}{{Partnership for policy integrity}}, \bibinfo{title}{Carbon emissions from burning biomass for energy}, \bibinfo{year}{2011}. \URLprefix \url{https://www.pfpi.net/wp-content/uploads/2011/04/PFPI-biomass-carbon-accounting-overview_April.pdf}.

\end{thebibliography}

\clearpage
\begin{appendices}

\clearpage
\section{Hierarchical clustering}
\label{ap: hierarchical clustering}

 \begin{table}[ht]
 \centering

\begin{tabular}{c|c|c|c}
        & Cluster 1        & Cluster 2       & Cluster 3       \\
Node    & weight: $\frac{145}{365}$  & weight: $\frac{117}{365}$ & weight: $\frac{103}{365}$ \\ \hline \hline
        & demand: day 43   & demand: day 263 & demand: day 177 \\
Finland & wind: day 12     & wind: day 257   & wind: day 204   \\
        & solar:   day 12  & solar: day 257  & solar: day 200  \\ \hline
        & demand: day 12   & demand: day 257 & demand: day 204 \\
Sweden  & wind: day 10     & wind: day 320   & wind: day 199   \\
        & solar:   day 6   & solar: day 91   & solar: day 159  \\ \hline
        & demand: day 12   & demand: day 257 & demand: day 204 \\
Norway  & wind: day 12     & wind: day 257   & wind: day 204   \\
        & solar:   day 6   & solar: day 77   & solar: day 210  \\ \hline
        & demand: day 12   & demand: day 257 & demand: day 204 \\
Denmark & wind: day 354    & wind: day 257   & wind: day 204   \\
        & solar: day 313   & solar: day 104  & solar: day 200  \\ \hline
        & demand: day 26   & demand: day 267 & demand: day 190 \\
Baltics & wind: day 88     & wind: day 258   & wind: day 141   \\
        & solar:   day 323 & solar: day 89   & solar: day 150  \\ \hline \hline
\end{tabular}
\caption{Representative days for each of the clusters and nodes}
  \label{Tab: Nordic_case_cluster}
\end{table}

\clearpage

\section{Alternative representations of the output data for Nordics case study}
\label{ap: alternative_representations}
\begin{figure}[h]
\centering
  \includegraphics[width=\linewidth]{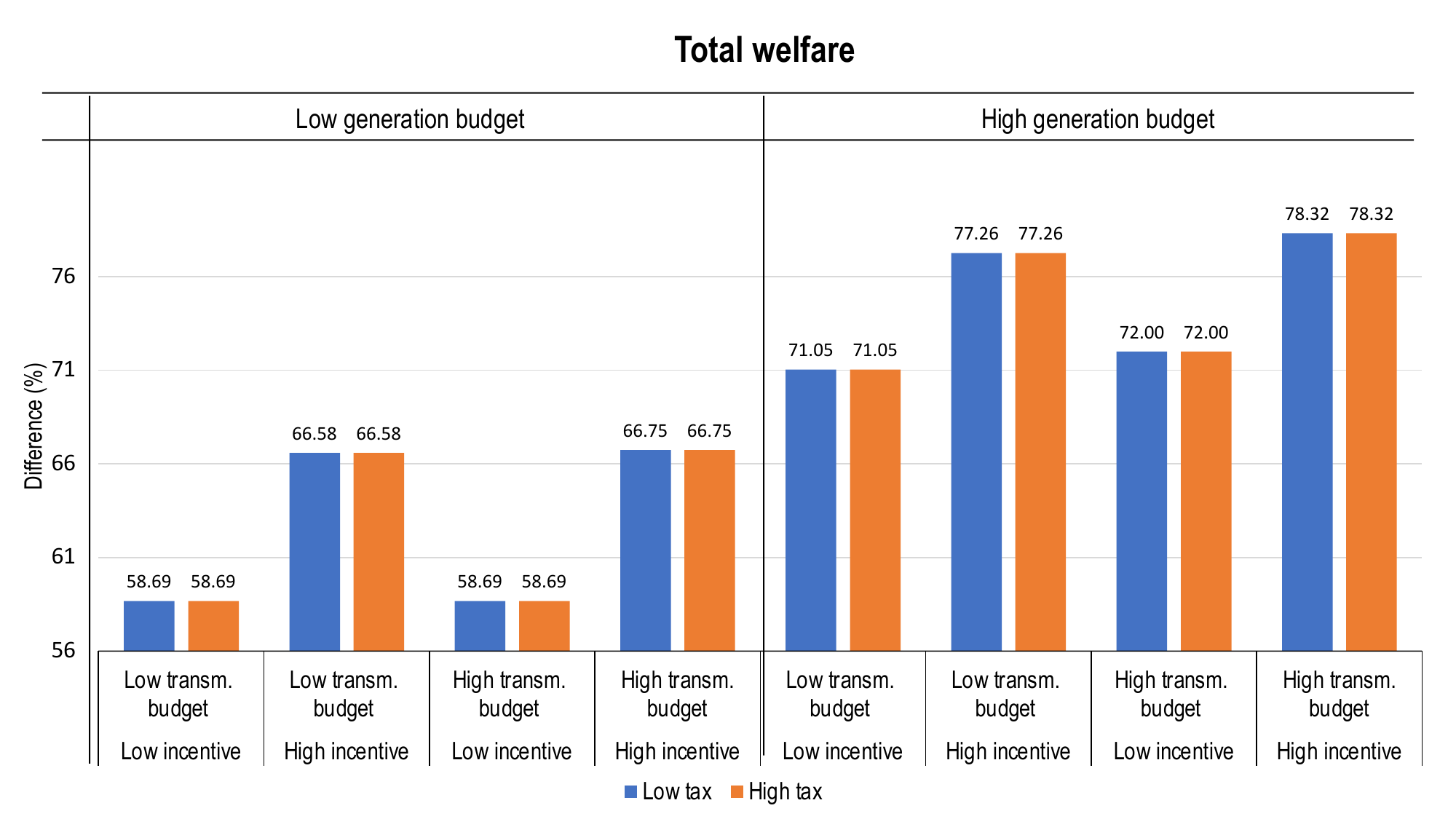}
    \caption{Case study: Impact of carbon taxes on welfare}
  \label{fig: nordics_tax_welfare}
 \end{figure}

\begin{figure}[h]
\centering
  \includegraphics[width=\linewidth]{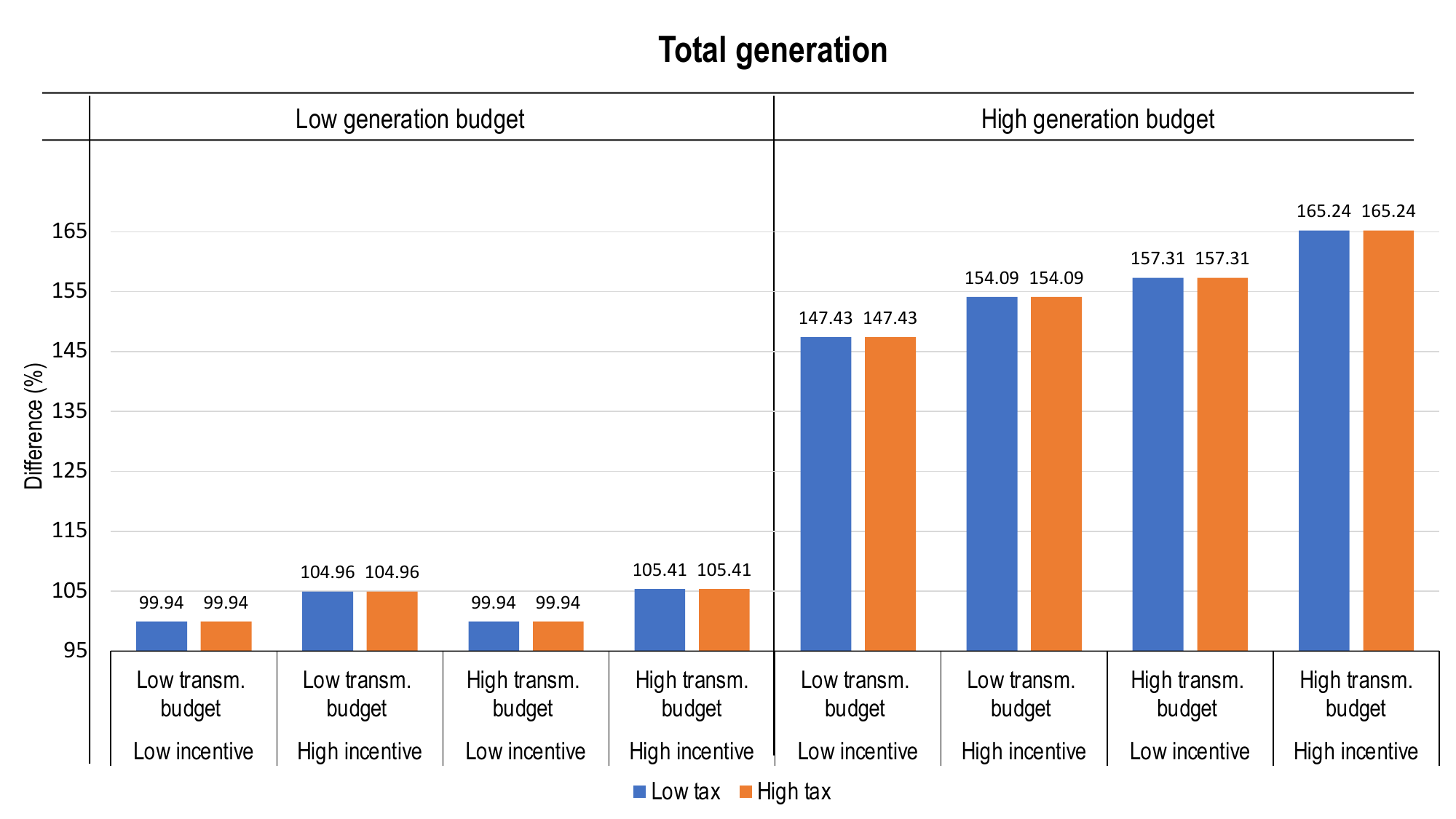}
   \caption{Case study: Impact of carbon taxes on total generation}
  \label{fig: nordics_tax_generation}
 \end{figure}

\begin{figure}[h]
\centering
  \includegraphics[width=\linewidth]{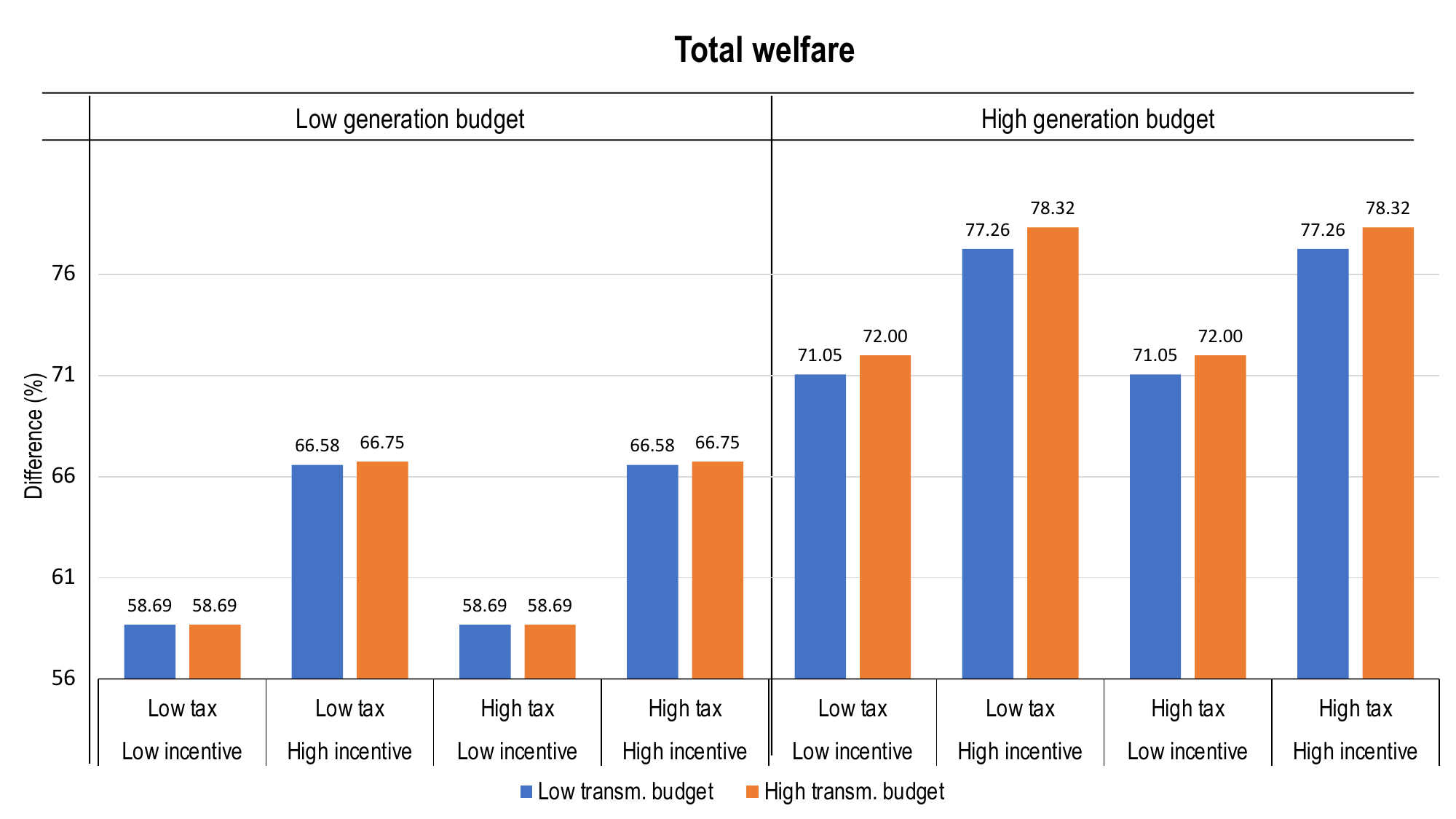}
    \caption{Case study: Impact of transmission capacity expansion budget on welfare}
  \label{fig: nordics_transm_welfare}
 \end{figure}

\begin{figure}[h]
\centering
  \includegraphics[width=\linewidth]{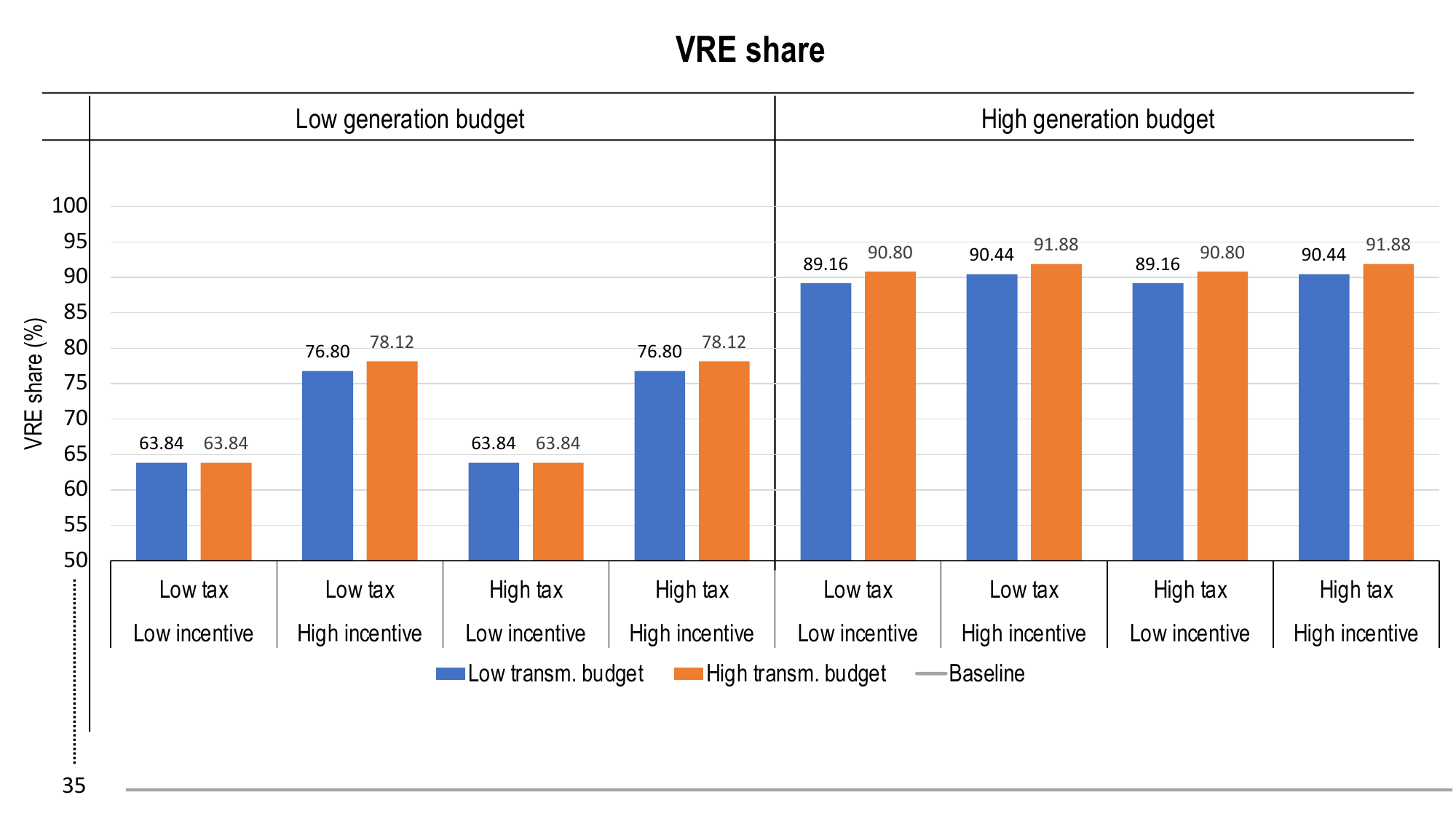}
  \caption{Case study: Impact of transmission capacity expansion budget on VRE share}
  \label{fig: nordics_transm_vres}
 \end{figure}

\begin{figure}[h]
\centering
  \includegraphics[width=\linewidth]{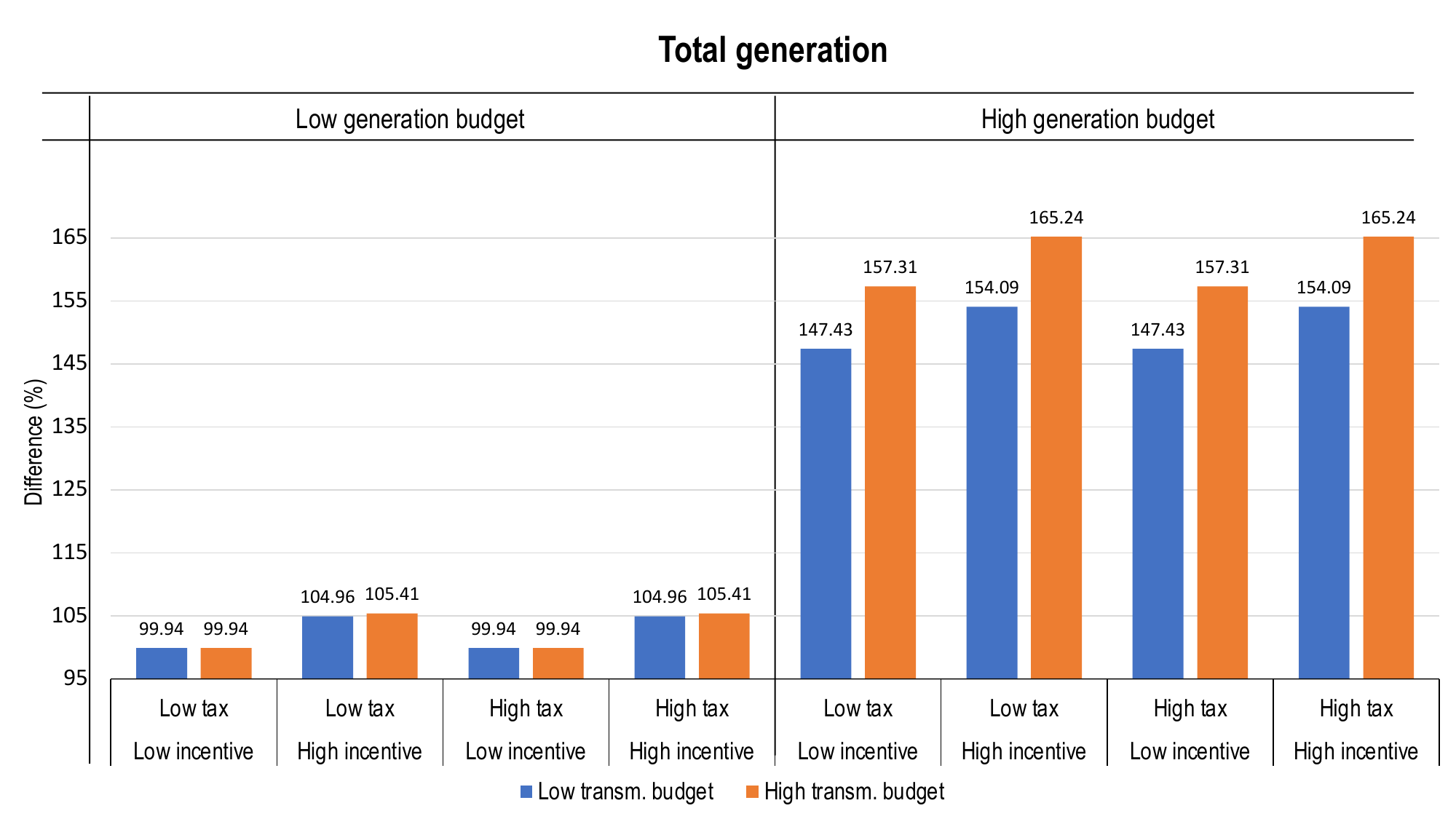}
   \caption{Case study: Impact of transmission capacity expansion budget on total generation}
  \label{fig: nordics_transm_generation}
 \end{figure}

 \begin{figure}[h]
\centering
  \includegraphics[width=\linewidth]{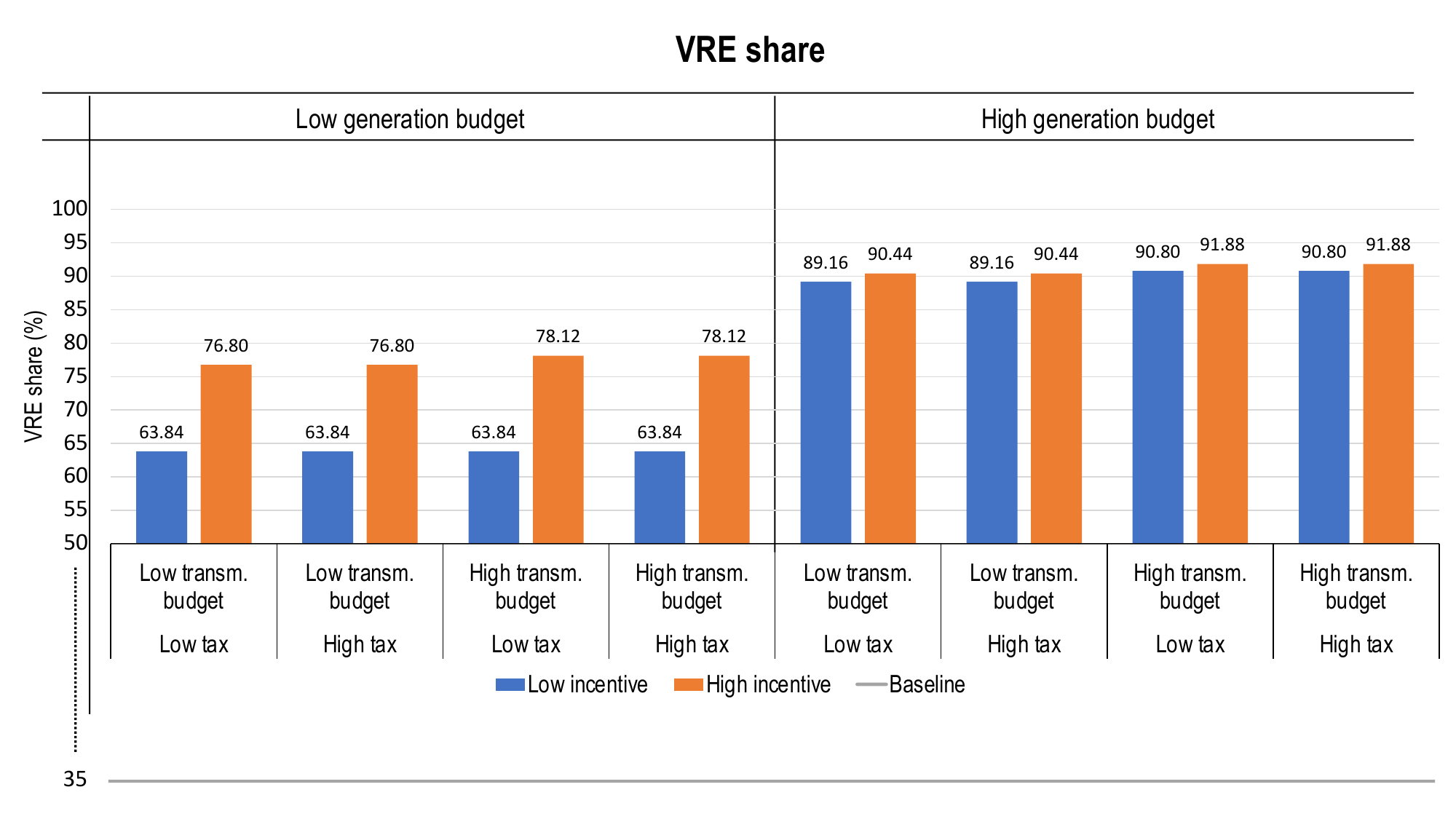}
  \caption{Case study: Impact of VRE investment incentive on VRE share}
  \label{fig: nordics_incentive_vres}
 \end{figure}
 
\clearpage

\section{List of Nomenclature}
\label{appendix: list of nomenclature}

% \begin{figure}[h]
% \centering
%   \includegraphics[width=\linewidth]{Pictures/general_parameters.pdf}
%     \caption{General parameters}
%  \end{figure}

\begin{table}[ht]
\centering
\begin{tabularx}{\linewidth}{c|L|s}
Symbol & Description                                                                                   & Units \\ \hline \hline
$T_t$   & Number of hours clustered for the time period $t \in T$                                                &    h   \\
$P_s$   & Probability of of the availability scenario   $s \in S$                                                &    \\
$D^{slp}_{s,n,t}$    & Slope of linear inverse demand function at scenario $s \in S$, node $n \in N$ in time period $t \in T$       &    \euro/MWh$^2$     \\
$D^{int}_{s,t,n}$    & Intercept of linear inverse demand function at scenario $s \in S$, node $n \in N$ in time period $t \in T$ &   \euro/MWh    \\
$K^{g+}_{i}$  & Generation capacity expansion investment budget for the producer $i \in I$                        &     \euro  \\  \hline \hline
\end{tabularx}
 \caption{General parameters}
\end{table}

% \begin{figure}[h]
% \centering
%   \includegraphics[width=\linewidth]{Pictures/primal_variables.pdf}
%     \caption{Primal variables}
%  \end{figure}

 \begin{table}[ht]
\centering
\begin{tabularx}{\linewidth}{c|L|s}
Symbol        & Description   & Units \\ \hline \hline
$g^e_{s,t,n,i}$ & Conventional generation of the type $e \in E$ at the node $n \in N$ by producer $i \in I$ considering scenario $s \in S$ and time period $t \in T$   &   MWh    \\
$g^r_{s,t,n,i}$ & VRE generation of the type $r \in R$ at the node $n \in N$ from producer $i \in I$ considering scenario $s \in S$ and time period $t \in T$ & MWh     \\
$q_{s,t,n}$ & Quantity of energy consumed at the node $n \in N$ considering scenario $s \in S$ during time period $t \in T$  &  MWh     \\
$f_{s,t,n,m}$ & Energy transferred from the node $n \in N$ to the node $m \in N$ considering scenario $s \in S$ and time period $t \in T$   &  MWh     \\
$l^+_{n,m}$ & Capacity added to the transmission line connecting nodes $n \in N$ and $m \in N$ & MW \\ 
$\overline{g}^{e+}_{n,i}$ & Generation capacity added to the conventional generation of the type $e \in E$ from the producer $i \in I$ at the node $n \in N$ & MW\\
$\overline{g}^{r+}_{n,i}$ & Generation capacity added to the VRE generation of the type $r \in R$ from the producer $i \in I$ at the node $n \in N$ &  MW \\ \hline \hline
\end{tabularx}
 \caption{Primal variables}
\end{table}

% \begin{figure}[h]
% \centering
%   \includegraphics[width=\linewidth]{Pictures/indices_and_sets.pdf}
%     \caption{Indices and sets}
%  \end{figure}

 \begin{table}[ht]
 \centering
\begin{tabular}{c|l}
Symbol & Description                       \\  \hline  \hline
$n \in N$    & Nodes                             \\
$s \in S$   & Availability scenarios            \\
$e \in E$    & Conventional energy sources       \\
$r \in R$    & Variable renewable energy sources \\
$i \in I$  & Power producer companies          \\
$t \in T$   & Time periods                      \\ \hline  \hline
\end{tabular}
 \caption{Indices and sets}
\end{table}

% \begin{figure}[h]
% \centering
%   \includegraphics[width=\linewidth]{Pictures/conventional_parameters.pdf}
%     \caption{Conventional generation parameters}
%  \end{figure}

\begin{table}[ht]
\centering
\begin{tabularx}{\linewidth}{c|L|s}
Symbol                & Description                                                                                                                              & Units \\ \hline \hline
$M^e_{n,i}$                   & Annualised maintenance cost for conventional generation of the type $e\in E$ from the producer $i \in I$ at the node $n \in N$                       &   \euro/MW    \\
$C^e_{n,i}$                  & Operational cost for conventional generation of the type $e\in E$  from producer  $i \in I$ the node $n \in N$                                       &  \euro/MWh     \\
$\overline{G}^e_{n,i}$                  & Installed generation capacity of the conventional generation of the type $e\in E$ from the producer $i \in I$ at the node $n \in N$                 & MW      \\
$I^e_{n,i}$                  & Annualised capacity expansion investment cost for the conventional generation of the type $e\in E$ from the producer $i \in I$  at the node $n \in N$  &  \euro/MW     \\
$R^{up, e}_{n,i}$              & Maximum ramp-up rate for the conventional generation of the type $e\in E$ at the node $n \in N$ from the producer $i \in I$                         &       \\
$R^{down, e}_{n,i}$ & Maximum ramp-down rate for the conventional generation of the type $e\in E$ at the node $n \in N$ from the producer $i \in I$                       &       \\
$D^e$ & Carbon tax for conventional generation of the type $e\in E$                                                                                  &  \euro/MWh     \\ \hline \hline
\end{tabularx}
  \caption{Conventional generation parameters}
\end{table}

% \begin{figure}[h]
% \centering
%   \includegraphics[width=\linewidth]{Pictures/vres_parameters.pdf}
%     \caption{VRE generation parameters}
%  \end{figure}

\begin{table}[ht]
\centering
\begin{tabularx}{\linewidth}{c|L|s}
Symbol        & Description   & Units \\ \hline \hline
$M^r_{n,i}$ & Annualised maintenance cost for VRE generation unit of the type $r \in  R$ from the producer $i \in I$ at the node $n \in N$ &   \euro/MW    \\
$\overline{G}^r_{n,i}$& Installed generation capacity of the VRE of the type $r \in  R$ from the producer $i \in I$  at the node $n \in N$  & MW      \\
$I^r_{n,i}$ & Annualised capacity expansion investment cost for the VRE of the type $r \in  R$ from the producer $i \in I$ at the node $n \in N$  &  \euro/MW     \\
$A^r_{s,t,n}$ & Availability factor for VRE type $r \in  R$ at the time period $t \in T$ considering scenario $s \in S$ at the node $n\in N$      &       \\ \hline \hline
\end{tabularx}
 \caption{VRE generation parameters}
\end{table}

% \begin{figure}[h]
% \centering
%   \includegraphics[width=\linewidth]{Pictures/dual_variables.pdf}
%     \caption{Dual variables}
%  \end{figure}

\begin{table}[ht]
\centering
\begin{tabularx}{\linewidth}{c|L|s}
Symbol        & Description   & Units \\ \hline \hline
$\theta_{s,t,n}$ & Shadow price on the power balance at node $n \in N$ considering scenario $s \in S$ during time period $t \in T$ &   \euro/MWh    \\
$\lambda^f_{s,t,n,m}$& Shadow price on the power flow primal feasibility constraint from the node $n \in N$ to the node $m \in N$ considering scenario $s \in S$ during time period $t \in T$ & \euro/MW      \\
$\beta^{f_1}_{s,t,n,m}$ & Shadow price on the transmission capacity for the power flow from the node $n \in N $ to the node $m \in N$ considering scenario $s \in S$ during time period $t \in T$ &  \euro/MW     \\
$\beta^{f_2}_{s,t,n,m}$ & Shadow price on the transmission capacity for the power flow from the node $n \in N $ to the node  $m \in N$  considering scenario $s \in S$ during time period $t \in T $   & \euro/MW  \\
$\beta^e_{s,t,n,i}$ & Shadow price on conventional energy capacity of the type $e \in E$ from the producer $i \in I$ at node $n \in N$ considering scenario $s \in S$ during time period $t \in T$  & \euro/MWh \\
$\beta^r_{s,t,n,i}$ & Shadow price on VRE capacity of the type $r \in R$ from the producer $i \in I$ at node $n \in N$ considering scenario $s \in S$ during time period $t \in T$  & \euro/MWh \\
$\beta^{up,e}_{s,t,n,i}$& Shadow price on the maximum ramp-up rate for the conventional generation of the type $e \in E$ at the node $n \in N$ from the producer $i \in I$ considering scenario $s \in S$ during time period $t \in T$  & \\
$\beta^{down,e}_{s,t,n,i}$& Shadow price on the maximum ramp-down rate for the conventional generation of the type $e \in E$ at the node $n \in N$ from the producer $i \in I$ considering scenario $s \in S$ during time period $t \in T$ & \\ 
$\lambda^{e+}_{n,i}$& Shadow price on non-negativity constraint for the conventional generation expansion of the type $e \in E$ at the node $n \in N$ from the producer $i \in I$  & \\
$\lambda^{r+}_{n,i}$& Shadow price on non-negativity constraint for the VRE generation expansion of the type $r \in R$ at the node $n \in N$ from the producer $i \in I$   & \\ 
$\lambda^e_{s,t,n,i}$& Shadow price on non-negativity constraint for the conventional generation of the type $e \in E$ at the node $n \in N$ from the producer $i \in I$ considering scenario $s \in S$ during time period $t \in T$  & \\
$\lambda^r_{s,t,n,i}$& Shadow price on non-negativity constraint for the VRE generation of the type $r \in R$ at the node $n \in N$ from the producer $i \in I$ considering scenario $s \in  S$ during time period $t \in T$  & \\
$\beta^{g+}_i$& Shadow price on upper-limit constraint for the generation expansion for the producer $i \in I$  & \\ \hline \hline
\end{tabularx}
 \caption{Dual variables}
\end{table}

% \begin{figure}[h]
% \centering
%   \includegraphics[width=\linewidth]{Pictures/transmission_parameters.pdf}
%     \caption{Transmission parameters}
%  \end{figure}

\begin{table}[ht]
\centering
\begin{tabularx}{\linewidth}{c|L|s}
Symbol        & Description   & Units \\ \hline \hline
$\overline{L}_{n,m}$ & Installed capacity at the line connecting nodes $n \in N$ and $m \in N$  &   MW    \\
$M^l_{n,m}$& Annualised maintenance cost for transmission per line connecting node $n \in N$ to the node $m \in N$   & \euro/MW      \\
$I^l_{n,m}$ & Annualised capacity expansion investment cost for the line connecting node $n \in N$ to the node $m \in N$    &  \euro/MW     \\
$K^{L+}$ & Capacity expansion investment budget   &  \euro     \\ \hline \hline
\end{tabularx}
 \caption{Transmission parameters}
\end{table}

\clearpage

\section{Numerical experiments: The data for the Nordic energy system case study}
\label{appendix: data}

To gather input data regarding demand and energy prices, we utilised the ENTSO-E Transparency Platform \citep{hirth_entso-e_2018}. Assuming each Nordic country was represented by one node, we had to aggregate price zones in some countries. For that, we accumulated the demand for all price zones in the country and calculated an average value for the energy price\footnote{We considered a simple moving average, i.e., the average price over the specified period. However, one could use other available alternatives, e.g., a volume-weighted average \citep{john1999technical}}. We performed the same aggregation for the Baltic countries.

We also used the ENTSO-E Transparency Platform to gather the data regarding the preinstalled generation capacity at each node. Analogously, we accumulated generation capacity by combining all price zones when necessary. Given the lack of information on the distribution of gas sources between open- and combined-cycle gas turbines, we assumed an equal share between them. Regarding the investment, variable operational, fixed maintenance and fuel costs, as well as the lifetime values for each of the energy sources but coal, we referred to \citep{KAN2020117015}. All the cost values for the coal are defined as in \citep{CHATZIMOURATIDIS2009778} while the lifetime was assumed to be 50 years in accordance with \citep{Cui_2019}. The values for the ramping levels regarding conventional energy sources are specified as in \citep{Rintamaki_2022}. The VRE availability factor represented as the share of the total installed capacity was defined following \citep{Condeixa_2022}.

We also accumulated all the data regarding hydropower generation capacity, including pumped storage, run-of-rivers and reservoirs provided by the ENTSO-E Transparency Platform. For the sake of simplicity, we assumed the hydropower capacity to be specified annually, implying that the hydropower reservoirs are refilled by the beginning of the next year. Another assumption involves the equal distribution of the total hydropower capacity among the time periods and scenarios. 

The transmission line investment costs, fixed maintenance costs and lifetime were specified as in \citep{KAN2020117015}. The installed transmission capacities were set according to \citep{NYCANDER2020942}. We accumulated transmission capacity for all connected price zones to define transmission line capacity between two nodes. Considering the cases when the transmission line allows for different voltages at the ends, to define a uniform capacity value, we considered an average value between them. Following \citep{SCHLACHTBERGER2017469}, the length of the transmission lines we considered to be the distances between geographical midpoints of each country. In the case of the Baltic countries, the closest Baltic country to the corresponding Nordic country was used to represent the Baltic countries node. However, when calculating the distance between Sweden and the Baltic countries node, the distance between the mid-point of Sweden and the mid-point of Lithuania was considered due to the existing transmission capacity between these two countries. For the same reason, the distance between the mid-point of Finland and the mid-point of Estonia was considered to represent the transmission line between Finland and the Baltic countries. 

The transmission expansion budget (TEB) values are based on the open data provided in the 2022 report of green solutions for Nordic energy system \citep{nordics_tsos}. The report states that the investments in the transmission infrastructure planned over the next ten years exceed \euro 25B. The project involves 4 TSOs: Fingrid (Finland), Energinet (Denmark), Statnett (Norway) and Svenska Kraftnät (Sweden). Hence, we assumed that a 10-year budget for one TSO would comprise \euro 6.25B. The generation expansion budget (GEB) values are based on the Swedish GenCo VATTENFALL report \citep{vattenfall_2022} indicating that the company planned growth investments of a total SEK 34 billion (about \euro3B) for 2022-2023. For the numerical experiments, the annual TEB and GEB values were levelled to a day-based scale by dividing them by 365. The data regarding the carbon tax is based on the data provided by the European review of energy taxation, carbon pricing and energy subsidies 2022 \citep{EU_review_taxation_2022}. Following \cite{ noauthor_GT, noauthor_coal, noauthor_biomass} we assumed the amount of emissions per MWh to be 0.50 tons of CO$_2$, 0.34 tons of CO$_2$, 0.99 tons of CO$_2$ and 0.23 tons of CO$_2$ for OCT, CCT, coal and biomass generation sources, respectively.

\end{appendices}

%\input{tables_figures}
% To print the credit authorship contribution details
%\printcredits

%% Loading bibliography style file
% \bibliographystyle{model1-num-names}

% Biography
%\bio{}
% Here goes the biography details.
%\endbio

%\bio{pic1}
% Here goes the biography details.
%\endbio

%% Authors are advised to use a BibTeX database file for their reference list.
%% The provided style file elsarticle-num.bst formats references in the required Procedia style

%% For references without a BibTeX database:

% \begin{thebibliography}{00}

%% \bibitem must have the following form:
%%   \bibitem{key}...
%%

% \bibitem{}

% \end{thebibliography}

\end{document}